\documentclass[a4paper]{article}
\usepackage{pdfsync}
\usepackage[utf8]{inputenc}
\usepackage{graphicx}
\usepackage{amssymb,amsfonts,amsmath,amsthm}
\usepackage{hyperref}
\usepackage{multirow}
\usepackage{verbatim}
\usepackage{color}
\usepackage{booktabs}
\usepackage[sectionbib,round]{natbib}
\usepackage{textcomp}
\usepackage{upquote}
\usepackage{bm}
\usepackage{dirtree}

\theoremstyle{plain}

\newtheorem{corollary}{Corollary}[section]
\newtheorem{lemma}{Lemma}[section]
\newtheorem{proposition}{Proposition}[section]
\theoremstyle{definition}
\newtheorem{definition}{Definition}[section]

\newtheorem{remark}{Remark}[section]

\newcommand*{\defeq}{\mathrel{\vcenter{\baselineskip0.5ex \lineskiplimit0pt \hbox{\scriptsize.}\hbox{\scriptsize.}}}=}
\newcommand{\minus}{\scalebox{0.6}{$-$}}
\newcommand{\plus}{\scalebox{0.6}{$+$}}

\begin{document}

\title{\textbf{A new family of models with generalized orientation in data envelopment analysis}}


\author{V. J. Bol\'os$^1$, R. Ben\'{\i}tez$^1$, V. Coll-Serrano$^2$ \\ \\
{\small $^1$ Dpto. Matem\'aticas para la Econom\'{\i}a y la Empresa, Facultad de Econom\'{\i}a.} \\
{\small $^2$ Dpto. Econom\'{\i}a Aplicada, Facultad de Econom\'{\i}a.} \\
{\small Universidad de Valencia. Avda. Tarongers s/n, 46022 Valencia, Spain.} \\
{\small e-mail\textup{: \texttt{vicente.bolos@uv.es}, \texttt{rabesua@uv.es}, \texttt{vicente.coll@uv.es}}} \\}

\date{June 2025}

\maketitle

\begin{abstract}
In the framework of data envelopment analysis, we review directional models \citep{Chambers1996,Chambers1998,Briec1997} and show that they are inadequate when inputs and outputs are improved simultaneously under constant returns to scale. Conversely, we introduce a new family of quadratically constrained models with generalized orientation and demonstrate that these models overcome this limitation. Furthermore, we extend the Farrell measure of technical efficiency using these new models. Additionally, we prove that the family of generalized oriented models satisfies some desired monotonicity properties. Finally, we show that the new models, although being quadratically constrained, can be solved through linear programs in a fundamental particular case.
\end{abstract}

\section{Introduction}

The concept of technical efficiency in production processes has its origins in the works of Koopmans \citep{Koop51}, Debreu \citep{Deb51} and Farrell \citep{Far57}. According to Debreu and Farrell, a measure of technical input efficiency is the minimum proportion by which a vector of inputs could be reduced while still producing a given output rate. In contrast, Koopmans considered an input combination to be efficient for a given output rate if a reduction in any input is not feasible for that output rate. In general, both concepts of technical efficiency are not coincident (Koopmans' concept is stronger), resulting in two ways of measuring efficiency.

In the framework of data envelopment analysis (DEA), the introduction of classic radial models \citep{Charnes1978,Charnes1979,Charnes1981} paved the way for the calculation of the Farrell measure of technical efficiency for a decision-making unit (DMU) that consumes inputs to produce outputs.
This measure indicates how much all inputs can be contracted (input oriented) or all outputs can be dilated (output oriented) while remaining within the production possibility set.
However, a limitation of the classic radial models is that input contractions or output dilations must be made in the same proportion.
To address this limitation, non-radial models and the so-called Russell measure were proposed in \cite{Fare1978}, following Koopmans' concept of technical efficiency and allowing for non-proportional improvements to reach the efficient frontier.
Subsequently, a number of models for estimating Russell-type efficiency measures have been proposed in the literature. Additive models were first introduced in \cite{Charnes1985}. Later, \cite{Pastor1999} and \cite{Tone2001} developed the slacks-based measure of efficiency (SBM) models, which exhibited the same characteristics as the additive models but provided an efficiency score.
Since then, non-radial models have been the focus of considerable attention from the scientific DEA community. Some recent works are \cite{Halicka2021, Gerami2022, Alves2023, Kuo2023, Bolos2023, Liu2025}, among others.
Nevertheless, all these non-radial models measure efficiency with respect to the ``farthest'' point in the portion of the efficient frontier that dominates the evaluated DMU. Consequently, they may be unsuitable for calculating a ``close'' efficient target. Some authors have adapted additive and SBM models to the search for the closest targets \citep{Apa2007, Tone2010}. However, the resulting models are not weakly monotonic in general, which gives rise to significant interpretation issues. Indeed, \cite{Fuk14} demonstrated that there is no weakly monotonic Russell-type score that employs a ratio-form least-distance approach to the closest projections over the efficient frontier.

Directional models were introduced in \cite{Chambers1996,Chambers1998,Briec1997} as an extension of radial models, allowing for non-proportional improvements of inputs and outputs while maintaining weak monotonicity. These models measure Farrell technical efficiency by customizing the proportions by which each variable is improved, taking account of particularities of the market and characterizing the criteria of management chosen by the producer. Some recent related  works are \cite{Sek2023, Pan2024, Bolos2024}. However, directional models are constructed in a linear way, and these improvements do not follow the constant returns to scale (CRS) assumption when inputs are contracted and outputs are dilated simultaneously (see Remark \ref{rem:crs}).

The objective of this paper is to define a novel family of \textit{generalized oriented models} that encompasses directional models and overcomes the aforementioned drawback under the CRS assumption. Generalized oriented models are capable of measuring Farrell technical efficiency in a manner analogous to directional models, with the ability to customize the proportions by which each variable is improved. 
Furthermore, the balance between input contractions and output dilations can be adjusted to be in accordance with the CRS assumption. However, these models can be quadratically constrained and hence, non-linear. In Section \ref{sec:mb}, we establish the mathematical background and in Section \ref{sec:gom}, we review directional models and introduce some novel generalized oriented models, called \textit{quadratic-CRS oriented models}. In Section \ref{sec:foe}, we utilize generalized oriented models to extend the Farrell measure of technical efficiency. In Section \ref{sec:mon}, we examine desirable monotonicity properties and in Section \ref{sec:zero}, we analyze the case of zeros in data. Although quadratic-CRS oriented models are non-linear, we demonstrate in Section \ref{sec:results} that they can be solved linearly in a fundamental particular case. Finally, we provide illustrative examples in Section \ref{sec:ex}, and present some conclusions in Section \ref{sec:fin}.

\section{Mathematical background}
\label{sec:mb}

Throughout this paper, we consider $\mathcal{D}=\left\{ \textrm{DMU}_1, \ldots ,\textrm{DMU}_n \right\} $ a set of $n$ DMUs with $m$ inputs and $s$ outputs. Matrices $X=(x_{ij})$ and $Y=(y_{rj})$ are the \emph{input} and \emph{output data matrices}, respectively, where $x_{ij}>0$ and $y_{rj}>0$ refers to the $i$-th input and $r$-th output, respectively, of the $j$-th DMU. Although we assume that the data are strictly positive, we study the case of zeros in the data set in Section \ref{sec:zero}.

In general, vector names are presented in boldface.
In this context, the zero vector and the all-ones vector are denoted by $\mathbf{0}$ and $\mathbf{1}$, respectively. The dimensions of these vectors are determined by the context. The term \emph{activity} is used to describe any pair of positive column vectors $\left( \mathbf{x},\mathbf{y}\right) $, where $\mathbf{x}\in \mathbb{R}^m_{>0}$ and $\mathbf{y}\in \mathbb{R}^s_{>0}$. Any $\textrm{DMU}_o \in \mathcal{D}$ has its associated activity $(\mathbf{x}_o,\mathbf{y}_o)$, where $\mathbf{x}_o=\left( x_{1o},\ldots ,x_{mo}\right) ^{\top}$ and $\mathbf{y}_o=\left( y_{1o},\ldots ,y_{so}\right) ^{\top}$ are positive column vectors. Therefore, a DMU can be identified with its activity in the same way that a point is identified with its coordinates.
Given two activities $\left( \mathbf{x},\mathbf{y} \right) $, $\left( \mathbf{x}',\mathbf{y}' \right) $, we say that $\left( \mathbf{x}',\mathbf{y}' \right) $ \emph{dominates} $\left( \mathbf{x},\mathbf{y} \right) $ if $\mathbf{x}'\leq \mathbf{x}$ and $\mathbf{y}\leq \mathbf{y}'$. The relation ``to be dominated by'' is a partial order on the set of activities.

The \emph{production possibility set} $P$ is the set of all the \emph{feasible activities} defined by $\mathcal{D}$. Assuming strong disposability of inputs and outputs, $P$ under CRS is the set formed by all the activities dominated by non-negative combinations of DMUs in $\mathcal{D}$:
\begin{equation}
\label{eq:p}
P =P(X,Y)\defeq \left\{ \left( \mathbf{x},\mathbf{y} \right) \in \mathbb{R}_{>0}^{m+s}\ \ | \ \ X\bm{\lambda}\leq \mathbf{x},\quad \mathbf{y}\leq Y\bm{\lambda},\quad \bm{\lambda}\geq \mathbf{0} \right\} ,
\end{equation}
where $\bm{\lambda}=(\lambda_1,\ldots,\lambda_n)^\top$ is a non-negative column vector. Other returns to scale can be considered by adding the following conditions in \eqref{eq:p}: $\sum _{j=1}^n \lambda _j=1$ for variable returns to scale (VRS), $0\leq \sum _{j=1}^n \lambda _j\leq 1$ for non-increasing returns to scale (NIRS), $\sum _{j=1}^n \lambda _j\geq 1$ for non-decreasing returns to scale (NDRS) or $L\leq \sum _{j=1}^n \lambda _j\leq U$ for generalized returns to scale (GRS), with $0\leq L\leq 1$ and $U\geq 1$.

We say that an activity $\left( \mathbf{x},\mathbf{y}\right) $ is \emph{efficient} (also known as \emph{strongly efficient} or \emph{Pareto-Koopmans efficient}) if there is not any other activity being feasible and dominating $\left( \mathbf{x},\mathbf{y}\right) $. 
The set of efficient activities in $P$ is called the \emph{efficient frontier} (also known as \emph{strongly efficient frontier} or \emph{Pareto-Koopmans frontier}) of $P$, and we denote it by $\partial ^{\text{S}}P$:
\begin{equation}
\label{eq:partialS}
\partial ^{\text{S}}P=\left\{ \left( \mathbf{x},\mathbf{y}\right) \in P\ \ | \ \ \left( \mathbf{x}',\mathbf{y}'\right) \in P,\ \mathbf{x}'\leq \mathbf{x},\ \mathbf{y}\leq \mathbf{y}'\ \Rightarrow \ \left( \mathbf{x}',\mathbf{y}'\right) = \left( \mathbf{x},\mathbf{y}\right) \right\}.
\end{equation}
We say that an activity $\left( \mathbf{x},\mathbf{y}\right) $ is \emph{weakly efficient} (also known as \emph{technically efficient} or \emph{Farrell efficient}) if there is not any other activity being feasible and improving all the variables of $\left( \mathbf{x},\mathbf{y}\right) $. The set of weakly efficient activities in $P$ is called the \emph{weakly efficient frontier} of $P$, and we denote it by $\partial ^{\text{W}}P$:
\begin{equation}
\label{eq:partialW}
\partial ^{\text{W}}P=\left\{ \left( \mathbf{x},\mathbf{y}\right) \in P\ \ | \ \ \mathbf{x}'<\mathbf{x},\ \mathbf{y}<\mathbf{y}' \ \Rightarrow \ \left( \mathbf{x}',\mathbf{y}'\right) \notin P \right\} ,
\end{equation}
where $<$ is a component-wise inequality, i.e. $x'_i<x_i$ for $i=1,\ldots ,m$, and $y_r<y'_r$ for $r=1,\ldots ,s$.

\begin{figure}[tb]
\begin{center}
  \includegraphics[width=0.3\textwidth]{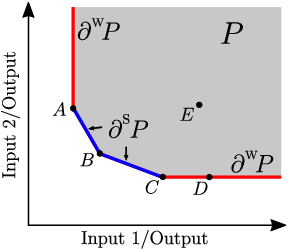}
\end{center}
\caption{Scheme of a set $\mathcal{D}=\left\{ A,B,C,D,E\right\}$ of 5 DMUs with 2 inputs and 1 output. The production possibility set $P$ is a closed set represented in grey (regardless of returns to scale), containing all the feasible activities defined by $\mathcal{D}$. The efficient frontier $\partial ^{\text{S}}P$ is represented by the blue line, and the weakly efficient frontier $\partial ^{\text{W}}P$ is obtained by the addition of the red lines. Efficient activities are those in $\partial ^{\text{S}}P$ but also those outside $P$. Efficient DMUs are $A,B,C$, while $E,D$ are not efficient. However, $D$ is weakly efficient.}
\label{fig0}
\end{figure}

It is clear that $\partial ^{\text{S}}P\subset \partial ^{\text{W}}P$, but they are not the same (see Fig. \ref{fig0}).
We have to remark that,  in this paper and in most of the literature, the term ``efficient'' refers to activities in $\partial ^{\text{S}}P$ that are efficient in the Koopmans' sense (i.e. strongly efficient), while ``weakly efficient'' refers to activities in $\partial ^{\text{W}}P$ that are technically efficient according to Debreu and Farrell, but not necessarily (strongly) efficient.
Moreover, we assume strong disposability of inputs and outputs and hence, the weakly efficient frontier of $P$ coincides with the boundary of $P$, i.e. $\partial ^{\text{W}}P=\partial P$ \citep{Fare1985}.
Note that this fact makes weakly efficient activities interesting because they compound the boundary of the production possibility set $P$, and they can be interpreted as the ``limits of the technology''.
Furthermore, assuming strictly positive data, we can apply a small improvement to any activity in $\partial ^{\text{W}}P$ to move it away from $P$ and transform it into an efficient one, although moving away from $P$ means that it is no longer feasible within the original technology.
For example, if we reduce the second input of DMU $D$ in Fig. \ref{fig0} by any small quantity, the resulting activity is efficient but it lies outside the original production possibility set $P$.

\section{Generalized oriented models}
\label{sec:gom}

Let $\textrm{DMU}_o=(\mathbf{x}_o,\mathbf{y}_o)\in \mathcal{D}$ be a DMU to evaluate. Directional models were originally defined in \cite{Chambers1996,Chambers1998} considering a \emph{direction} of improvement $\left( -\mathbf{g}^{\minus},\mathbf{g}^{\plus}\right) \neq \mathbf{0}$ with $\mathbf{g}^{\minus}\in \mathbb{R}^m$ and $\mathbf{g}^{\plus}\in \mathbb{R}^s$ being non-negative. These models generalize classic radial models measuring the maximum displacement from $(\mathbf{x}_o,\mathbf{y}_o)$ along the direction $\left( -\mathbf{g}^{\minus},\mathbf{g}^{\plus}\right)$ remaining inside $P$. In other words, directional models tell us how much the variables of $\textrm{DMU}_o$ can be improved following the preassigned direction and still be feasible within the original technology. If $\left( -\mathbf{g}^{\minus},\mathbf{g}^{\plus}\right) =\left( -\mathbf{x}_o,\mathbf{0}\right)$, we obtain classic \emph{input oriented radial} models; while if $\left( -\mathbf{g}^{\minus},\mathbf{g}^{\plus}\right) =\left( \mathbf{0},\mathbf{y}_o\right)$, we obtain classic \emph{output oriented radial} models \citep{Charnes1978}.

An alternative formulation was introduced in \cite{Briec1997} using a generalized \emph{orientation} of the form $\left( \mathbf{d}^{\minus},\mathbf{d}^{\plus}\right) \neq \mathbf{0}$ instead of $\left( -\mathbf{g}^{\minus},\mathbf{g}^{\plus}\right) $, with the \emph{orientation vectors} $\mathbf{d}^{\minus}\in \mathbb{R}^m$ and $\mathbf{d}^{\plus}\in \mathbb{R}^s$ being non-negative. The relation between the orientation $\left( \mathbf{d}^{\minus},\mathbf{d}^{\plus}\right) $ and the direction $\left( -\mathbf{g}^{\minus},\mathbf{g}^{\plus}\right) $ of the original directional models is the following: $g^{\minus}_i=d^{\minus}_i x_{io}$ and $g^{\plus}_r=d^{\plus}_r y_{ro}$ for $i=1,\ldots ,m$ and $r=1,\ldots ,s$.

In this paper, we are going to define a new family of models which includes directional models making use of a generalized orientation $\left( \mathbf{d}^{\minus},\mathbf{d}^{\plus}\right) $. In all these \emph{generalized oriented models}, the orientation coefficients characterize the criteria of management chosen by the producer, representing the relative ``ease of improvement'' of each variable and thus, determining an improvement strategy: the higher the orientation coefficient of a variable, the easier it is to improve that variable; if the orientation coefficient of a variable is equal to $0$, then that variable can not be changed. In this way, $\mathbf{d}^{\minus}$ determines the proportions by which inputs are contracted and $\mathbf{d}^{\plus}$ determines the proportions by which outputs are dilated, but the balance between input contractions and output dilations will be determined by the model (see Remark \ref{rem:no}).

All the generalized oriented models that we are going to define calculate a \emph{target} activity  $(\mathbf{x}_o^*,\mathbf{y}_o^*)\in \partial ^{\text{W}}P$ indicating the maximum improvement of  $(\mathbf{x}_o,\mathbf{y}_o)$ following the strategy determined by the orientation while remaining feasible. Since the target dominates $\textrm{DMU}_o$, it can be written as
\begin{equation}
\label{eq:projslack}
\left( \mathbf{x}_o^*,\mathbf{y}_o^*\right) =\left( \mathbf{x}_o-\mathbf{t}^{\minus *},\mathbf{y}_o+\mathbf{t}^{\plus *}\right) ,
\end{equation}
where the non-negative column vectors $\mathbf{t}^{\minus *}\in \mathbb{R}^m$ and $\mathbf{t}^{\plus *}\in \mathbb{R}^s$ are the \emph{input} and \emph{output target slacks vectors} of $\textrm{DMU}_o$, respectively. We define the \emph{relative input} and \emph{output target slacks vectors} of $\textrm{DMU}_o$, denoted by $\bm{\tau}^{\minus *}$ and $\bm{\tau}^{\plus *}$ respectively, as the column vectors with coefficients
\begin{equation}
\label{eq:taudef}
\tau ^{\minus *}_i\defeq t^{\minus *}_i/x_{io},\quad \tau ^{\plus *}_r\defeq t^{\plus *}_r/y_{ro},\qquad i=1,\ldots ,m,\quad r=1,\ldots ,s.
\end{equation}
Note that, from \eqref{eq:taudef}, it holds $0\leq \tau ^{\minus *}_i<1$ and $\tau ^{\plus *}_r\geq 0$. On the other hand, the target \eqref{eq:projslack} can also be written as
\begin{equation}
\label{eq:projcondil}
\left( \mathbf{x}_o^*,\mathbf{y}_o^*\right) =\left( \mathrm{diag}\left( \bm{\theta}^*\right) \mathbf{x}_o,\mathrm{diag}\left( \bm{\phi}^*\right) \mathbf{y}_o\right) ,
\end{equation}
where $\mathrm{diag}\left( \bm{\theta}^*\right) $ and $\mathrm{diag}\left( \bm{\phi}^*\right) $ are diagonal contraction and dilation matrices, respectively, with
\begin{equation}
\label{eq:thetagen}
\bm{\theta}^*=1-\bm{\tau}^{\minus *},\qquad \bm{\phi}^*=1+\bm{\tau}^{\plus *},
\end{equation}
being the \emph{input target contractions} and \emph{output target dilations vectors} of $\textrm{DMU}_o$, respectively. Note that, from \eqref{eq:thetagen}, it holds $0<\theta ^*_i\leq 1$ and $\phi ^*_r \geq 1$, for $i=1,\ldots ,m$ and $r=1,\ldots ,s$.

\begin{remark}[\textbf{Input and output oriented radial models}]
Since the coefficients of the orientation $(\mathbf{d}^{\minus},\mathbf{d}^{\plus})$ represent the \textit{relative} ease of improvement of each variable, proportional orientations must give the same target. 
Taking this into account, we say that a generalized oriented model is \emph{input oriented radial} if $(\mathbf{d}^{\minus},\mathbf{d}^{\plus})$ is proportional to $(\mathbf{1},\mathbf{0})$. In this case, we must obtain the same targets as in the classic input oriented radial models. On the other hand, we say that a generalized oriented model is \emph{output oriented radial} if $(\mathbf{d}^{\minus},\mathbf{d}^{\plus})$ is proportional to $(\mathbf{0},\mathbf{1})$. In this case, we must obtain the same targets as in the classic output oriented radial models.
\end{remark}

\begin{remark}[\textbf{Non-oriented radial models. Balance}]
\label{rem:no}
In the particular case where $(\mathbf{d}^{\minus},\mathbf{d}^{\plus})$ is proportional to $(\mathbf{1},\mathbf{1})$, we say that a generalized oriented model is \emph{non-oriented radial}. In this scenario, there is not any discriminated variable in the sense that all variables have the same technical difficulty to be improved. As a consequence, all inputs must be contracted with the same contraction coefficient $0<\theta ^*\leq 1$, and all outputs must be dilated with the same dilation coefficient $\phi ^*\geq 1$. The relationship between $\theta ^*$ and $\phi ^*$ determines the model's own \emph{balance} between input contractions and output dilations in the calculation of the target.
\end{remark}

\begin{remark}[\textbf{Balanced according to CRS}]
\label{rem:crs}
Under the CRS assumption, if an activity multiplies its inputs by $\lambda >0$, then it can produce $\lambda $ times its original outputs. We define an equivalence relation for activities based on this assumption: $\left( \mathbf{x},\mathbf{y}\right) \sim \left( \mathbf{x}',\mathbf{y}'\right) $ if and only if $\exists \lambda >0$ such that $\left( \mathbf{x},\mathbf{y}\right) =\lambda \left( \mathbf{x}',\mathbf{y}'\right) $.
The equivalence class of an activity $\left( \mathbf{x},\mathbf{y}\right) $ are all the feasible activities that can be generated from $\left( \mathbf{x},\mathbf{y}\right) $, assuming CRS. In other words, two activities are equivalent if and only if they individually generate the same production possibility set \eqref{eq:p} under CRS. Hence, given $0<\theta \leq 1$, we have
$\left( \theta \mathbf{x},\mathbf{y}\right) \sim \frac{1}{\theta}\left( \theta \mathbf{x},\mathbf{y}\right) =\left( \mathbf{x},\frac{1}{\theta}\mathbf{y}\right) $.
That is, contracting all the inputs with a contraction coefficient $\theta $ is equivalent to dilating all the outputs with a dilation coefficient $1/\theta $. This fact is reflected in the classic input and output oriented radial models under CRS, for which the output projection dilation $\phi ^*$ of a DMU evaluated by the output oriented radial model is the inverse of the input projection contraction $\theta ^*$ of that DMU evaluated by the input oriented radial model.
Taking this into account, we say that a generalized oriented model is \emph{balanced according to CRS} if it holds
$\phi ^*=\frac{1}{\theta ^*}$
in the non-oriented radial case (i.e. $(\mathbf{d}^{\minus},\mathbf{d}^{\plus})$ proportional to $(\mathbf{1},\mathbf{1})$), where $\theta ^*$ is the input target contraction coefficient (applied to all inputs) and $\phi ^*$ is the output target dilation coefficient (applied to all outputs). 
\end{remark}

\begin{remark}[\textbf{Efficient projection and second stage}]
\label{rem:ss}
Once applied a generalized oriented model to $\textrm{DMU}_o=(\mathbf{x}_o,\mathbf{y}_o)\in \mathcal{D}$ and calculated the target $(\mathbf{x}^*_o,\mathbf{y}^*_o)$, this weakly efficient activity may not be efficient (i.e. in $\partial ^{\text{S}}P$). In this case, the model also computes an \emph{efficient projection} as an activity in $\partial ^{\text{S}}P$ dominating $(\mathbf{x}_o^*,\mathbf{y}_o^*)$, and the vectors resulting from the difference between the efficient projection and the target are called  \emph{inefficiency slacks vectors}. However, if there is more than one activity in $\partial ^{\text{S}}P$ dominating the target, a \emph{second stage} can be performed in order to establish a criterion for selecting the efficient projection. This second stage usually consists on computing the so called \emph{max-slack solution} by means of a classic additive model applied to $(\mathbf{x}^*_o,\mathbf{y}^*_o)$.
\end{remark}

\subsection{Linear oriented models}

In this paper, we will refer to the original directional models defined in \cite{Chambers1996,Chambers1998,Briec1997} as \emph{linear oriented} (LO) models. Given $\textrm{DMU}_o=(\mathbf{x}_o,\mathbf{y}_o)\in \mathcal{D}$ and an orientation $\left( \mathbf{d}^{\minus},\mathbf{d}^{\plus}\right) $, LO models calculate the target \eqref{eq:projcondil} so that $\bm{\theta}^*$ and $\bm{\phi}^*$ are optimal for the program
\begin{equation}
\label{eq:dirgenlin}
\def\arraystretch{1.2}
\begin{array}[t]{rl}
\beta ^*_{\text{L}}=\max \limits_{\beta ,\bm{\theta},\bm{\phi}} & \beta \\
\text{s.t.} & \left( \mathrm{diag}\left( \bm{\theta}\right) \mathbf{x}_o,\mathrm{diag}\left( \bm{\phi}\right) \mathbf{y}_o\right) \in P, \\
& \bm{\theta}=1-\beta\mathbf{d}^{\minus}, \\
& \bm{\phi}=1+\beta\mathbf{d}^{\plus}, \\
& 
\beta \geq 0.
\end{array}
\end{equation}
Note that the target \eqref{eq:projcondil} is, effectively, in $ \partial ^{\text{W}}P$. The optimal value $\beta ^*_{\text{L}}$ satisfies $0\leq \beta ^*_{\text{L}}< 1/\left\| \mathbf{d}^{\minus}\right\| _{\infty}$, except for the case of output orientation ($\mathbf{d}^{\minus}=\mathbf{0}$), in which $\beta ^*_{\text{L}}$ is not upper bounded. Moreover, $\beta ^*_{\text{L}}=0$ if and only if $\textrm{DMU}_o$ is its own target.
In the particular case where $(\mathbf{d}^{\minus},\mathbf{d}^{\plus})$ is proportional to $(\mathbf{1},\mathbf{0})$, we obtain the same targets as in classic input oriented radial models; on the other hand, if $(\mathbf{d}^{\minus},\mathbf{d}^{\plus})$ is proportional to $(\mathbf{0},\mathbf{1})$, we obtain the same targets as in classic output oriented radial models.

The input target contractions and output target dilations vectors are given by
\begin{equation}
\label{eq:thetalin}
\bm{\theta}^* = 1-\beta ^*_{\text{L}}\mathbf{d}^{\minus},
\qquad \bm{\phi}^* = 1+\beta ^*_{\text{L}}\mathbf{d}^{\plus},
\end{equation}
and the relative target slacks vectors are given by
\begin{equation}
\label{eq:taulin}
\bm{\tau}^{\minus *} = \beta ^*_{\text{L}}\mathbf{d}^{\minus},
\qquad \bm{\tau}^{\plus *} = \beta ^*_{\text{L}}\mathbf{d}^{\plus}.
\end{equation}
It follows from \eqref{eq:taulin} that the relative target slacks vectors are directly proportional to the orientation vectors, with the same proportionality factor $\beta ^*_{\text{L}}$ (see also Fig. \ref{fig1}). On the other hand, LO models are not balanced according to CRS (see Remark \ref{rem:crs}) because, in the non-oriented case, the output target dilation coefficient $1+\beta ^*_{\text{L}}$ is not the inverse of the input target contraction coefficient $1-\beta ^*_{\text{L}}$, except for the trivial case $\beta ^*_{\text{L}}=0$ (see \eqref{eq:thetalin}).

Taking into account \eqref{eq:p}, program \eqref{eq:dirgenlin} under CRS can be expressed as a linear program:
\begin{equation}
\label{eq:dirmodlin}
\def\arraystretch{1.2}
\begin{array}[t]{rll}
\beta ^*_{\text{L}}=\max \limits_{\beta ,\bm{\lambda}} & \beta \\
\text{s.t.} & \sum_{j=1}^n\lambda _jx_{ij}+\beta d^{\minus}_ix_{io}\leq x_{io},& \quad i=1,\ldots ,m, \\
& \sum_{j=1}^n\lambda _jy_{rj}-\beta d^{\plus}_ry_{ro}\geq y_{ro},& \quad r=1,\ldots ,s, \\
& \beta \geq 0,\quad \bm{\lambda}\geq \mathbf{0}.
\end{array}
\end{equation}
The efficient projection is given by $(X\bm{\lambda}^*,Y\bm{\lambda}^*)\in \partial ^{\text{S}}P$, and the inefficiency slacks vectors are given by $\mathbf{s}^{\minus *} =\mathrm{diag}\left( \bm{\theta}^*\right) \mathbf{x}_o-X\bm{\lambda}^*$ and $\mathbf{s}^{\plus *} =Y\bm{\lambda}^*-\mathrm{diag}\left( \bm{\phi}^*\right) \mathbf{y}_o$, where $\bm{\lambda}^*$ is optimal for program \eqref{eq:dirmodlin}. If there is not uniqueness for $\bm{\lambda}^*$, a second stage can be performed (see Remark \ref{rem:ss}).
Finally, \eqref{eq:dirmodlin} can be adapted for different returns to scale by adding the corresponding constraints: $\sum _{j=1}^n \lambda _j=1$ (VRS), $0\leq \sum _{j=1}^n \lambda _j\leq 1$ (NIRS), $\sum _{j=1}^n \lambda _j\geq 1$ (NDRS) or $L\leq \sum _{j=1}^n \lambda _j\leq U$ (GRS).

\subsection{Quadratic-CRS oriented models}

Given $\textrm{DMU}_o=(\mathbf{x}_o,\mathbf{y}_o)\in \mathcal{D}$ and an orientation $\left( \mathbf{d}^{\minus},\mathbf{d}^{\plus}\right) $, \emph{quadratic-CRS oriented} (QO) models calculate the target \eqref{eq:projcondil} so that $\bm{\theta}^*$ and $\bm{\phi}^*$ are optimal for the program
\begin{equation}
\label{eq:dirgenquad}
\def\arraystretch{1.2}
\begin{array}[t]{rl}
\beta ^*_{\text{Q}}=\max \limits_{\beta ,\bm{\theta},\bm{\phi}} & \beta \\
\text{s.t.} & \left( \mathrm{diag}\left( \bm{\theta}\right) \mathbf{x}_o,\mathrm{diag}\left( \bm{\phi}\right) \mathbf{y}_o\right) \in P, \\
& \bm{\theta}=1-\beta \mathbf{d}^{\minus}, \\
& \bm{\phi}=1/\left( 1-\beta \mathbf{d}^{\plus}\right) , \\
& 
\beta \geq 0,
\end{array}
\end{equation}
where the quotient of a vector is component-wise. Note that the target \eqref{eq:projcondil} is, effectively, in $ \partial ^{\text{W}}P$. The optimal value $\beta ^*_{\text{Q}}$ satisfies $0\leq \beta ^*_{\text{Q}}< 1/\left\| (\mathbf{d}^{\minus},\mathbf{d}^{\plus}) \right\| _{\infty}$.
Moreover, $\beta ^*_{\text{Q}}=0$ if and only if $\textrm{DMU}_o$ is its own target, as in LO models.
In the particular case where $(\mathbf{d}^{\minus},\mathbf{d}^{\plus})$ is proportional to $(\mathbf{1},\mathbf{0})$, we obtain the same targets as in classic input oriented radial models; on the other hand, if $(\mathbf{d}^{\minus},\mathbf{d}^{\plus})$ is proportional to $(\mathbf{0},\mathbf{1})$, we obtain the same targets as in classic output oriented radial models.

\begin{figure}[tb]
\begin{center}
  \includegraphics[width=0.65\textwidth]{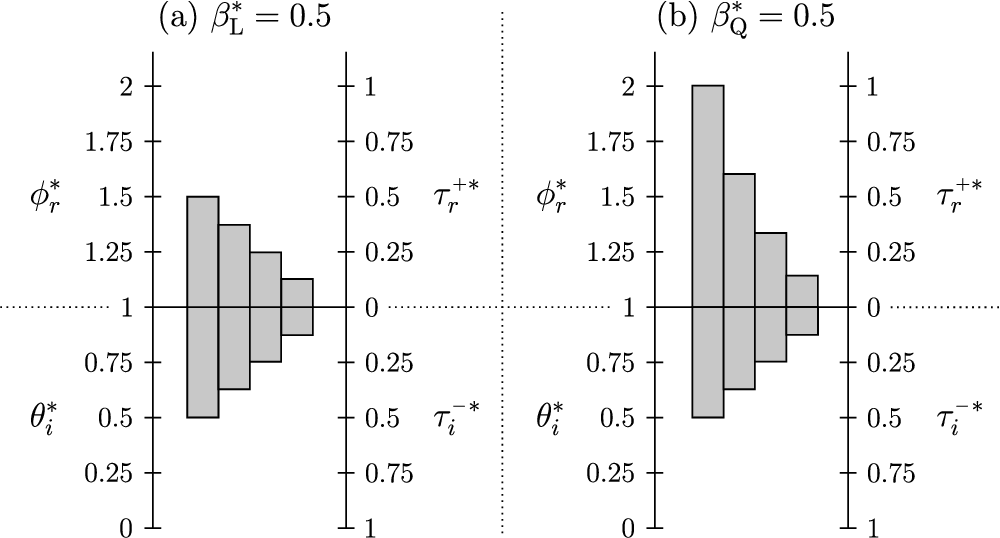}
\end{center}
\caption{Input target contractions and output target dilations, made by (a) LO and (b) QO models. Considering DMUs with $m=4$ inputs and $s=4$ outputs, and orientation vectors $\mathbf{d}^{\protect \minus}=\mathbf{d}^{\protect \plus}=(1,0.75,0.5,0.25)$, we have represented on the left axis input target contractions, $\theta ^*_i$ for $i=1,\ldots ,4$, and output target dilations, $\phi ^*_r$ for $r=1,\ldots ,4$, corresponding to optimal solutions with (a) $\beta ^*_{\text{L}}=0.5$, and (b) $\beta ^*_{\text{Q}}=0.5$, according to \eqref{eq:thetalin} and \eqref{eq:thetaquad}, respectively. Moreover, according to \eqref{eq:thetagen}, the length of a bar is the value of the corresponding relative target slack, represented on the right axis.}
\label{fig1}
\end{figure}

\begin{figure}[tb]
\begin{center}
  \includegraphics[width=0.65\textwidth]{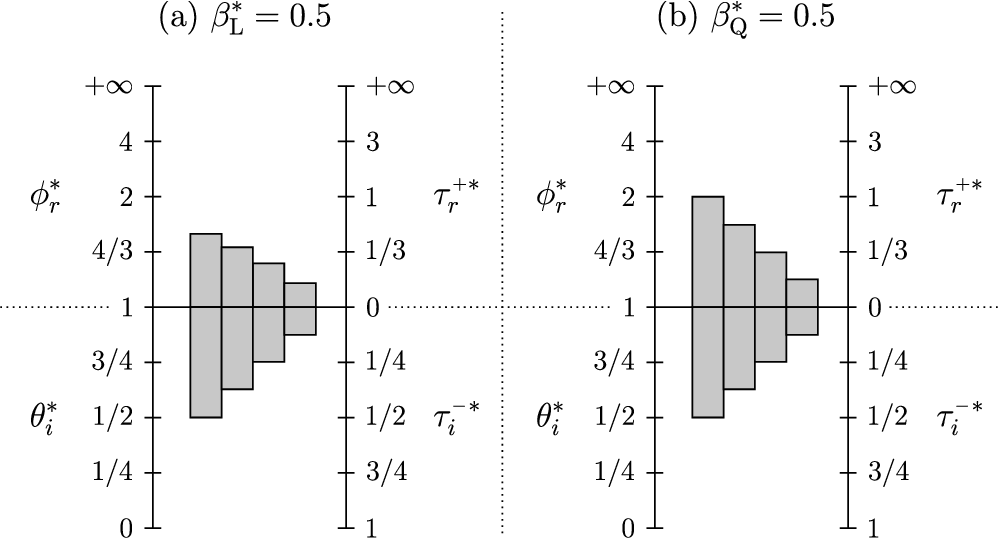}
\end{center}
\caption{The same plots represented in Figure \ref{fig1}, but the scale of dilations is inverse with respect to the scale of contractions, which continues being linear. It is noted that QO models are according to CRS because the target dilation of an output with orientation coefficient $d$ is the inverse of the target contraction of an input with the same orientation coefficient $d$. This property can be seen in (b), where the dilation bars are an exact reflection of the contraction bars.}
\label{fig2}
\end{figure}

The input target contractions and output target dilations vectors are given by
\begin{equation}
\label{eq:thetaquad}
\bm{\theta}^* = 1-\beta ^*_{\text{Q}}\mathbf{d}^{\minus},
\qquad \bm{\phi}^* = 1/\left( 1-\beta ^*_{\text{Q}}\mathbf{d}^{\plus}\right) ,
\end{equation}
and the relative projection slacks vectors are given by
\begin{equation}
\label{eq:tauquad}
\bm{\tau}^{\minus *} = \beta ^*_{\text{Q}}\mathbf{d}^{\minus},
\qquad \bm{\tau}^{\plus *} = \beta ^*_{\text{Q}}\mathbf{d}^{\plus}/\left( 1-\beta ^*_{\text{Q}}\mathbf{d}^{\plus}\right) ,
\end{equation}
where the quotient of vectors is component-wise. It follows from \eqref{eq:tauquad} that $\bm{\tau}^{\plus *}$ and  $\mathbf{d}^{\plus}$ are not proportional, unlike LO models (see Fig. \ref{fig1}). On the other hand, QO models are balanced according to CRS (see Remark \ref{rem:crs}) because the target dilation of an output with orientation coefficient $d$ is the inverse of the target contraction of an input with the same orientation coefficient $d$ (see Fig. \ref{fig2}).

\begin{proposition}
\label{propab}
Let us consider the optimal values $\beta ^*_{\text{Q}}$ and $\beta ^*_{\text{L}}$ from QO and LO models given by \eqref{eq:dirgenquad} and \eqref{eq:dirgenlin}, respectively. Then, $\beta ^*_{\text{Q}} \leq \beta ^*_{\text{L}}$.
Moreover, $\beta ^*_{\text{Q}}=\beta ^*_{\text{L}}$ for any evaluated DMU if and only if $\mathbf{d}^{\plus}=\mathbf{0}$.
\end{proposition}

Taking into account \eqref{eq:p} and $0\leq \beta ^*_{\text{Q}}<1$, program \eqref{eq:dirgenquad} under CRS can be written as a quadratically constrained linear program with a convex feasible set:
\begin{equation}
\label{eq:dirmodquad}
\def\arraystretch{1.2}
\begin{array}[t]{rll}
\beta ^*_{\text{Q}}=\max \limits_{\beta ,\bm{\phi},\bm{\lambda}} & \beta \\
\text{s.t.} & \sum_{j=1}^n\lambda _jx_{ij}+\beta d^{\minus}_ix_{io}\leq x_{io},& \quad i=1,\ldots ,m, \\
& \sum_{j=1}^n\lambda _jy_{rj}-\phi _r y_{ro}\geq 0,& \quad r=1,\ldots ,s, \\
& \left( \beta d^{\plus}_r-1\right) \phi _r\leq -1,& \quad r=1,\ldots ,s, \\
& \beta \geq 0,\quad \bm{\phi}\geq \mathbf{0},\quad  \bm{\lambda}\geq \mathbf{0}.
\end{array}
\end{equation}
Note that the constraint $\left( \beta d^{\plus}_r-1\right) \phi _r\leq -1$ is equivalent to $\left( \beta d^{\plus}_r-1\right) \phi _r=-1$ because we want to maximize $\beta $. Moreover, although the functions $\left( \beta d^{\plus}_r-1\right) \phi _r$ are not convex, the constraints $\left( \beta d^{\plus}_r-1\right) \phi _r\leq -1$ jointly with $\bm{\phi}\geq \mathbf{0}$ define a convex set.
The efficient projection is given by $(X\bm{\lambda}^*,Y\bm{\lambda}^*)\in \partial ^{\text{S}}P$ and the inefficiency slacks vectors are given by $\mathbf{s}^{\minus *} =\mathrm{diag}\left( \bm{\theta}^*\right) \mathbf{x}_o-X\bm{\lambda}^*$ and $\mathbf{s}^{\plus *} =Y\bm{\lambda}^*-\mathrm{diag}\left( \bm{\phi}^*\right) \mathbf{y}_o$, where $\bm{\lambda}^*$ is optimal for program \eqref{eq:dirmodquad}. If there is not uniqueness for $\bm{\lambda}^*$, a second stage can be performed (see Remark \ref{rem:ss}).
Finally, although it is designed for CRS, program \eqref{eq:dirmodquad} can be adapted for different returns to scale by adding the corresponding constraints: $\sum _{j=1}^n \lambda _j=1$ (VRS), $0\leq \sum _{j=1}^n \lambda _j\leq 1$ (NIRS), $\sum _{j=1}^n \lambda _j\geq 1$ (NDRS) or $L\leq \sum _{j=1}^n \lambda _j\leq U$ (GRS).

\subsection{Orientation and improvement cost function}

In this section, we are going to differentiate between \emph{controllable} variables, that can be improved at some cost, and \emph{uncontrollable} variables, that can not be improved or the user do not want to improve them.
Moreover, we suppose that there exists a hypothetical improvement cost function of the form $f\left( \bm{\tau}^{\minus}_{\textrm{c}},\bm{\tau}^{\plus}_{\textrm{c}}\right) $, that gives the cost of applying the positive relative slacks vectors $\left( \bm{\tau}^{\minus}_{\textrm{c}},\bm{\tau}^{\plus} _{\textrm{c}}\right)$ to the controllable variables of the evaluated $\textrm{DMU}_o$. We are going to suppose that $f$ is unknown, but we know its gradient $\nabla f\left( \mathbf{0},\mathbf{0}\right) $, with all its components being strictly positive.

Precisely, the components of $\nabla f\left( \mathbf{0},\mathbf{0}\right) $ quantify the marginal cost of improving each controllable variable of $\textrm{DMU}_o$. Hence, we can take any vector proportional to $1/\nabla f(\mathbf{0},\mathbf{0})$ (where the inverse of the gradient is component-wise) as the orientation for controllable variables, $\left( \mathbf{d}^{\minus}_{\textrm{c}},\mathbf{d}^{\plus}_{\textrm{c}}\right) $. On the other hand, we take $\left( \mathbf{d}^{\minus}_{\textrm{u}},\mathbf{d}^{\plus}_{\textrm{u}}\right) =\left( \mathbf{0},\mathbf{0}\right) $ for uncontrollable variables. In this case, each orientation coefficient effectively can be interpreted as the relative ease of improvement of the corresponding variable according to the improvement cost function.

The next proposition gives us an interpretation of the parameters $\beta ^*_{\mathrm{L}}$ and $\beta ^*_{\mathrm{Q}}$ in terms of the cost of transforming the evaluated $\textrm{DMU}_o$ into its target, and suggesting that an appropriate choice of the orientation is precisely \eqref{eq:imprLQO}.

\begin{proposition}
\label{prop:impr}
Let $f\left( \bm{\tau}^{\minus}_{\textrm{c}},\bm{\tau}^{\plus}_{\textrm{c}}\right) $ be an improvement cost function with respect to controllable variables, and let us take
\begin{equation}
\label{eq:imprLQO}
\left( \mathbf{d}^{\minus}_{\textrm{c}},\mathbf{d}^{\plus}_{\textrm{c}}\right) =\frac{1}{m_{\textrm{c}}+s_{\textrm{c}}}\frac{1}{\nabla f\left( \mathbf{0},\mathbf{0}\right) },
\end{equation}
for controllable variables, where $m_{\textrm{c}}+s_{\textrm{c}}$ is the total number of controllable variables, and the inverse of the gradient is component-wise. On the other hand, we take $\left( \mathbf{d}^{\minus}_{\textrm{u}},\mathbf{d}^{\plus}_{\textrm{u}}\right) =\left( \mathbf{0},\mathbf{0}\right) $ for uncontrollable variables. Then, $\beta ^*_{\mathrm{L}}$ (from LO model \eqref{eq:dirgenlin}) and $\beta ^*_{\mathrm{Q}}$ (from QO model \eqref{eq:dirgenquad}) are linear approximations of the cost of transforming $\textrm{DMU}_o$ into its corresponding target.
\end{proposition}

Finally, if we want orientation coefficients between $0$ and $1$, we must take $\left( \mathbf{d}^{\minus}_{\textrm{c}},\mathbf{d}^{\plus}_{\textrm{c}}\right) =\| \nabla f\left( \mathbf{0},\mathbf{0}\right) \|_{\infty} /\nabla f\left( \mathbf{0},\mathbf{0}\right) $. In this case, we will obtain the same target as in Proposition \ref{prop:impr}, but we have to multiply $\beta ^*_{\mathrm{L}}$ or $\beta ^*_{\mathrm{Q}}$ by $\| \nabla f\left( \mathbf{0},\mathbf{0}\right) \|_{\infty} \cdot \left( m_{\textrm{c}}+s_{\textrm{c}}\right) $ in order to obtain a linear approximation of the cost of transforming $\textrm{DMU}_o$ into its target.

\section{Farrell oriented efficiency}
\label{sec:foe}

The \emph{Farrell measure of technical efficiency} (or \emph{Farrell efficiency}, for short) is an efficiency score developed from the concept of technical efficiency introduced by Debreu and Farrell \citep{Deb51,Far57}. It measures efficiency with respect to the weak efficient frontier $\partial ^{\text{W}}P$, as distinguished from the so-called \emph{Russell-type efficiencies} which measure efficiency with respect to the strongly efficient frontier $\partial ^{\text{S}}P$, following the Koopmans' concept of efficiency \citep{Koop51}.

Classic radial models were constructed in \cite{Charnes1978} for computing the Farrell efficiency in DEA. In the input oriented case, the Farrell efficiency is given by the input target contraction, $\theta ^*$ (where $0<\theta ^*\leq 1$ such that $\bm{\theta}^*=\left( \theta ^*,\ldots ,\theta ^*\right) $); while in the output oriented case, it is given by the inverse of the output target dilation, $1/\phi ^*$ (where $\phi ^*\geq 1$ such that  $\bm{\phi}^*=\left( \phi ^*,\ldots ,\phi ^*\right) $).

With the aim of extending the radial Farrell efficiency to generalized oriented models, we are going to define an efficiency score inspired by the \emph{enhanced Russell graph (ERG) measure of efficiency}  \citep{Pastor1999}, also known as the \emph{slacks-based measure (SBM) of efficiency} \citep{Tone2001}.

\begin{definition}
\label{def:fardir}
Given $\textrm{DMU}_o=(\mathbf{x}_o,\mathbf{y}_o)\in \mathcal{D}$ and an orientation $\left( \mathbf{d}^{\minus},\mathbf{d}^{\plus}\right) $, the \emph{Farrell oriented efficiency} of $\textrm{DMU}_o$ with orientation $\left( \mathbf{d}^{\minus},\mathbf{d}^{\plus}\right) $ is defined as
\begin{equation}
\label{eq:rho1}
\rho ^*\defeq \frac{\frac{1}{m}\sum _{i=1}^m\theta ^*_i}{\frac{1}{s}\sum _{r=1}^s\phi ^*_r},
\end{equation}
where $\bm{\theta}^*$ and $\bm{\phi}^*$ are the input target contractions and output target dilations vectors of $\textrm{DMU}_o$, respectively, calculated from a generalized oriented model with orientation $\left( \mathbf{d}^{\minus},\mathbf{d}^{\plus}\right) $. Taking into account \eqref{eq:thetagen}, an equivalent expression is given by
\begin{equation}
\label{eq:rho2}
\rho ^*=\frac{1-\frac{1}{m}\sum _{i=1}^m \tau ^{\minus *}_i}{1+ \frac{1}{s}\sum _{r=1}^s \tau ^{\plus *}_r},
\end{equation}
where $\bm{\tau}^{\minus *},\bm{\tau}^{\plus *}$ are the corresponding relative target slacks vectors.
\end{definition}

We have that $0<\rho ^*\leq 1$, and $\rho ^*=1$ if and only if $\bm{\tau}^{\minus *}=\mathbf{0}$ and $\bm{\tau}^{\plus *}=\mathbf{0}$, i.e. if and only if $\textrm{DMU}_o$ is its own target. Moreover, the Farrell oriented efficiency \eqref{eq:rho1} clearly generalizes the input and output oriented Farrell efficiencies of classic radial models: $\rho ^*=\theta ^*$ if $\bm{\theta}^*=\left( \theta ^*,\ldots ,\theta ^*\right) $, $\bm{\phi}^*=\left( 1,\ldots ,1\right) $ (input oriented radial case), and $\rho ^*=1/\phi ^*$ if $\bm{\theta}^*=\left( 1,\ldots ,1\right) $, $\bm{\phi}^*=\left( \phi ^*,\ldots ,\phi ^*\right) $ (output oriented radial case).
In fact, \eqref{eq:rho1} is the ratio between the average contraction of inputs and the average dilation of outputs.

Given an orientation $\left( \mathbf{d}^{\minus},\mathbf{d}^{\plus}\right) $, the corresponding Farrell oriented efficiency of $\left( \mathbf{x}_o,\mathbf{y}_o\right) $ can be computed for LO and QO models. Let us denote these scores by $\rho ^*_{\text{L}}$ and $\rho ^*_{\text{Q}}$, respectively. For LO models, taking into account \eqref{eq:thetalin}, $\rho ^*_{\text{L}}$ can be written as
\begin{equation}
\label{eq:rholin}
\rho ^*_{\text{L}}\left( \beta ^*_{\text{L}};\mathbf{d}^{\minus};\mathbf{d}^{\plus}\right) =\displaystyle \frac{1-\beta ^*_{\text{L}}\frac{1}{m}\sum _{i=1}^m d^{\minus}_i}{1+\beta ^*_{\text{L}}\frac{1}{s}\sum _{r=1}^s d^{\plus}_r},
\end{equation}
where $\beta ^*_{\text{L}}$ is determined by \eqref{eq:dirgenlin}.
As particular cases of \eqref{eq:rholin}, we have
\begin{itemize}
\item
$\rho ^*_{\text{L}}\left( \beta ^*_{\text{L}};1,\ldots ,1;1,\ldots ,1\right) =\displaystyle \frac{1-\beta ^*_{\text{L}}}{1+\beta ^*_{\text{L}}}$ for non-oriented radial models,
\vspace{0.2cm}
\item
$\rho ^*_{\text{L}}\left( \beta ^*_{\text{L}};1,\ldots ,1;0,\ldots ,0\right) =1-\beta ^*_{\text{L}}$ for input oriented radial models,
\vspace{0.2cm}
\item
$\rho ^*_{\text{L}}\left( \beta ^*_{\text{L}};0,\ldots ,0;1,\ldots ,1\right) =\displaystyle \frac{1}{1+\beta ^*_{\text{L}}}$ for output oriented radial models.
\end{itemize}
For QO models, taking into account \eqref{eq:thetaquad}, $\rho ^*_{\text{Q}}$ can be written as
\begin{equation}
\label{eq:rhoquad}
\rho ^*_{\text{Q}}\left( \beta ^*_{\text{Q}};\mathbf{d}^{\minus};\mathbf{d}^{\plus}\right) =\displaystyle \frac{1-\beta ^*_{\text{Q}}\frac{1}{m}\sum _{i=1}^m d^{\minus}_i}{\displaystyle \frac{1}{s}\sum _{r=1}^s \frac{1}{1-\beta ^*_{\text{Q}}d^{\plus}_r}},
\end{equation}
where $\beta ^*_{\text{Q}}$ is determined by \eqref{eq:dirgenquad}. As particular cases of \eqref{eq:rhoquad}, we have
\begin{itemize}
\item
$\rho ^*_{\text{Q}}\left( \beta ^*_{\text{Q}};1,\ldots ,1;1,\ldots ,1\right) =\left( 1-\beta ^*_{\text{Q}}\right) ^2$ for non-oriented radial models,
\vspace{0.2cm}
\item
$\rho ^*_{\text{Q}}\left( \beta ^*_{\text{Q}};1,\ldots ,1;0,\ldots ,0\right) =\rho ^*_{\text{Q}}\left( \beta ^*_{\text{Q}};0,\ldots ,0;1,\ldots ,1\right) =1-\beta ^*_{\text{Q}}$ for input and output oriented radial models.
\end{itemize}

As shown in Fig. \ref{fig1} and \ref{fig2}, LO models are not balanced according to CRS (see Remark \ref{rem:crs}). This leads to the fact that the expression \eqref{eq:rholin} of $\rho ^*_{\text{L}}$ is different for input oriented radial models and for output oriented radial models.
On the other hand, the expression of $\rho ^*_{\text{Q}}$ for input oriented radial models coincides with the expression of $\rho ^*_{\text{Q}}$ for output oriented radial models.

\section{Monotonicity}
\label{sec:mon}

The relation ``to be dominated by'' is a partial order on the set of activities. In this section, we are going to show that generalized oriented models are monotonic with respect to this partial order, as well as the Farrell oriented efficiencies. It should be noted that we are going to study non-strict monotonicity, also known as weak monotonicity, since strict (or strong) monotonicity is not fulfilled by classic radial models and Farrell efficiency.

Let us consider $\mathcal{D}=\left\{ \textrm{DMU}_1, \ldots ,\textrm{DMU}_n \right\} $ a set of \emph{reference DMUs}, with $X,Y$ being the input and output data matrices, respectively, and defining a production possibility set $P$ given by \eqref{eq:p}. We are going to use generalized oriented models to define functions on $P$: given $(\mathbf{x},\mathbf{y})\in P$ and an orientation $(\mathbf{d}^{\minus},\mathbf{d}^{\plus})$, we define $\beta ^*(\mathbf{x},\mathbf{y})$ (with respect to $\mathcal{D}$ and with orientation $(\mathbf{d}^{\minus},\mathbf{d}^{\plus})$) as the value of $\beta ^*$ that would return a generalized oriented model evaluating a new hypothetical DMU with activity $(\mathbf{x},\mathbf{y})$, i.e. considering $\mathcal{D}\cup \left\{ (\mathbf{x},\mathbf{y})\right\} $ as the set of DMUs. For LO models, taking into account \eqref{eq:dirmodlin}, we have
\begin{equation}
\label{eq:dirmodlinxy}
\def\arraystretch{1.2}
\begin{array}[t]{rll}
\beta ^*_{\text{L}}(\mathbf{x},\mathbf{y})=\max \limits_{\beta ,\bm{\lambda}} & \beta \\
\text{s.t.} & \sum_{j=1}^n\lambda _jx_{ij}-(1-\beta d^{\minus}_i)x_i\leq 0,& \quad i=1,\ldots ,m, \\
& \sum_{j=1}^n\lambda _jy_{rj}-(1+\beta d^{\plus}_r)y_r\geq 0,& \quad r=1,\ldots ,s, \\
& \beta \geq 0,\quad \bm{\lambda}\geq \mathbf{0}.
\end{array}
\end{equation}
Note that, since $(\mathbf{x},\mathbf{y})\in P$, there exist activities in $P$ of the form $(X\bm{\lambda},Y\bm{\lambda})$ dominating $(\mathbf{x},\mathbf{y})$. Hence, \eqref{eq:dirmodlinxy} is well defined and $\beta ^*_{\text{L}}(\mathbf{x},\mathbf{y})$ is, effectively, the LO model with orientation $(\mathbf{d}^{\minus},\mathbf{d}^{\plus})$ evaluating a DMU with activity $(\mathbf{x},\mathbf{y})$. On the other hand, program \eqref{eq:dirmodlinxy} is unfeasible for $(\mathbf{x},\mathbf{y})\notin P$, but in this case $(\mathbf{x},\mathbf{y})$ is efficient and $\beta ^*_{\text{L}}(\mathbf{x},\mathbf{y})$ could be taken equal to $0$. Analogously, we can define $\beta ^*_{\text{Q}}(\mathbf{x},\mathbf{y})$ for  $(\mathbf{x},\mathbf{y})\in P$ considering \eqref{eq:dirmodquad}:
\begin{equation}
\label{eq:dirmodquadxy}
\def\arraystretch{1.2}
\begin{array}[t]{rll}
\beta ^*_{\text{Q}}(\mathbf{x},\mathbf{y})=\max \limits_{\beta ,\bm{\lambda}} & \beta \\
\text{s.t.} & \sum_{j=1}^n\lambda _jx_{ij}-(1-\beta d^{\minus}_i)x_i\leq 0,& \quad i=1,\ldots ,m, \\
\vspace{-0.4cm} \\
& \sum_{j=1}^n\lambda _jy_{rj}-\left( \displaystyle \frac{1}{1-\beta d^{\plus}_r}\right) y_r\geq 0,& \quad r=1,\ldots ,s, \\
\vspace{-0.5cm} \\
& \beta \geq 0,\quad \bm{\lambda}\geq \mathbf{0}.
\end{array}
\end{equation}
Functions \eqref{eq:dirmodlinxy} and \eqref{eq:dirmodquadxy} are defined under CRS, but they can be adapted to different returns to scale in the usual way and the results are not affected. Moreover, \eqref{eq:dirmodquadxy} can be expressed as a quadratically constrained linear program with convex feasible set.

\begin{proposition}
\label{prop:monbeta}
Let $\beta ^*(\mathbf{x},\mathbf{y})$ be a function on $P$ defined by $\beta ^*_{\text{L}}(\mathbf{x},\mathbf{y})$ or $\beta ^*_{\text{Q}}(\mathbf{x},\mathbf{y})$. Then,
$\frac{\partial \beta ^*}{\partial x_i}(\mathbf{x},\mathbf{y})\geq 0$ and $\frac{\partial \beta ^*}{\partial y_r}(\mathbf{x},\mathbf{y})\leq 0$ for $i=1,\ldots ,m$ and $r=1,\ldots ,s$. 
\end{proposition}

Proposition \ref{prop:monbeta} assures that the value of the parameter $\beta ^*$ returned by a generalized oriented model evaluating a DMU does not increase if the activity of the DMU is improved. In other words, generalized oriented models are monotonically decreasing considering the partial order ``to be dominated by''.

Likewise, we can also define functions on $P$ by means of Farrell oriented efficiencies: given $(\mathbf{x},\mathbf{y})\in P$ and an orientation $(\mathbf{d}^{\minus},\mathbf{d}^{\plus})$, we define $\rho ^*(\mathbf{x},\mathbf{y})$ (with respect to $\mathcal{D}$ and orientation $(\mathbf{d}^{\minus},\mathbf{d}^{\plus})$) as
\begin{equation}
\label{eq:rhoxy}
\rho ^*(\mathbf{x},\mathbf{y})=\rho ^*\left( \beta ^*(\mathbf{x},\mathbf{y});\mathbf{d}^{\minus};\mathbf{d}^{\plus}\right).
\end{equation}
Note that we can take $\bm{\rho}^*(\mathbf{x},\mathbf{y})=1$ for $(\mathbf{x},\mathbf{y})\notin P$ and hence, $\bm{\rho}^*$ is continuous in the set of all activities.

\begin{proposition}
\label{prop:monrho}
Let $\rho ^*(\mathbf{x},\mathbf{y})$ be a function on $P$ defined by \eqref{eq:rhoxy} with $\beta ^*(\mathbf{x},\mathbf{y})$ given by \eqref{eq:dirmodlinxy} or \eqref{eq:dirmodquadxy}. Then,
$\frac{\partial \rho ^*}{\partial x_i}(\mathbf{x},\mathbf{y})\leq 0$ and $\frac{\partial \rho ^*}{\partial y_r}(\mathbf{x},\mathbf{y})\geq 0$ for $i=1,\ldots ,m$ and $r=1,\ldots ,s$. 
\end{proposition}

Proposition \ref{prop:monrho} assures that the Farrell oriented efficiency of a DMU does not decrease if the activity of the DMU is improved. In other words, Farrell oriented efficiencies are monotonically increasing considering the partial order ``to be dominated by''.

\section{Zeros in data}
\label{sec:zero}

We have assumed in Section \ref{sec:mb} that the input and output data matrices $X,Y$ are strictly positive. However, programs \eqref{eq:dirmodlin} and \eqref{eq:dirmodquad} can deal with zeros in data, both in the DMU under evaluation ($\text{DMU}_o$) and in the other DMUs of $\mathcal{D}$. Nevertheless, the orientation $\left( \mathbf{d}^{\minus},\mathbf{d}^{\plus}\right) $ must have, at least, one non-zero coefficient corresponding to a non-zero variable of $\text{DMU}_o$, since the orientation coefficient of a zero variable becomes irrelevant in those programs.

For the case of zeros in input data, suppose that the first input of $\text{DMU}_o$ is zero, i.e. $x_{1o}=0$. Then, the first input of the target must be also zero and hence, $s^{\minus *}_1=0$. The relative target slack $\tau ^{\minus *}_1$ is taken equal to $0$ for convenience and so, the target contraction coefficient $\theta ^*_1$ is equal to $1$.
With respect to the Farrell oriented efficiency, the parameter $m$ (number of inputs) in the expression of $\rho ^*$, \eqref{eq:rho1} or \eqref{eq:rho2}, must be reduced to the number of non-zero inputs \citep{Tone2001}.

On the other hand, for the case of zeros in output data, suppose that the first output of $\text{DMU}_o$ is zero, i.e. $y_{1o}=0$.
Following \cite{Tone2001}, we differentiate between two cases:

\textit{Case 1.} $\text{DMU}_o$ has no possibility to produce the first output. Then, $s^{\plus *}_1$ and $\tau ^{\minus *}_1$ are taken equal to $0$ for convenience and so, the target dilation coefficient $\phi ^*_1$ is taken equal to $1$. With respect to the Farrell oriented efficiency, the parameter $s$ (number of outputs) in the expression of $\rho ^*$, \eqref{eq:rho1} or \eqref{eq:rho2}, must be reduced by $1$.

\textit{Case 2.} $\text{DMU}_o$ has the potential to produce the first output. Then, we may replace $y_{1o}$ by a sufficiently small positive number, e.g.
$\min\left\{ y_{1j}\,\, |\,\, y_{1j}>0,\,\, j=1,\ldots ,m\right\}/10$.
In this case, the term $\tau ^{\plus *}_1$ (whose magnitude is controlled by $d^{\plus}_1$) in \eqref{eq:rho2} has the role of a penalty.

With respect to negative data, we have to take into account that the assumption of CRS is not possible in technologies under negative data. Anyway, we have different options about dealing with negative data. The simplest one is to consider negative inputs as positive outputs and negative outputs as positive inputs \citep{Scheel2001}. Another method is to make all negative data positive by adding a big enough scalar to the negative variables \citep{Charnes1983}. Nevertheless, this approach can be only applied if the model has translation invariance, and generalized oriented models are not invariant under translations. Finally, \citet{Alla2015} proposed linear directional models that, under non-constant returns to scale, can deal with negative data. More recently, \citet{Tavana2021} expanded these models considering flexible measures.

\section{Results from particular cases}
\label{sec:results}

Generalized oriented models collapse into classic radial models if $(\mathbf{d}^{\minus},\mathbf{d}^{\plus})$ is proportional to $(\mathbf{1},\mathbf{0})$ (input oriented radial), or $(\mathbf{d}^{\minus},\mathbf{d}^{\plus})$ is proportional to $(\mathbf{0},\mathbf{1})$ (output oriented radial). Furthermore, if $\mathbf{d}^{\plus}=\mathbf{0}$, all generalized oriented models become the same.
Note that, with respect to the relation between LO and QO models, we can only affirm that there exists an orientation vector $\mathbf{d}^{\plus}_{\text{L}}$ for which the target calculated by a QO model with orientation $( \mathbf{0},\mathbf{d}^{\plus}_{\text{Q}}) $ coincides with the target calculated by the corresponding LO model with orientation $( \mathbf{0},\mathbf{d}^{\plus}_{\text{L}}) $. However, $\mathbf{d}^{\plus}_{\text{L}}$ can not be computed \textit{a priori}, i.e. before computing $\beta ^*_{\text{Q}}$.
Notwithstanding this, there are some other important particular cases in which we can solve quadratically constrained models by means of linear programs, as we are going to see next.

In the rest of this section, we are going to study the particular case $d^{\minus}_i=d^{\minus}$ for $i=1,\ldots ,m$, and $d^{\plus}_r=d^{\plus}$ for $r=1,\ldots ,s$, with $d^{\minus},d^{\plus}> 0$, under CRS. In this fundamental case, there is no discrimination between inputs nor between outputs in the sense that all inputs have the same orientation coefficient $d^{\minus}$ and all outputs have the same orientation coefficient $d^{\plus}$, but these two coefficients can be different. Since these coefficients represent the \textit{relative} ease of improvement of inputs and outputs, we could assume that they are normalized by  $\| \, \| _{\infty}$ and hence $\max \left\{d^{\minus},d^{\plus}\right\} =1$; however, all the following results do not require this assumption.
Note that we are excluding $d^{\minus}=0$ or $d^{\plus}=0$ because, in this case, generalized oriented models coincide with the classic radial models. On the other hand, if $d^{\minus}=d^{\plus}=1$ we have the important case of non-oriented radial models, which is included in our study.

In the particular case $d^{\minus}_i=d^{\minus}$ for $i=1,\ldots ,m$, and $d^{\plus}_r=d^{\plus}$ for $r=1,\ldots ,s$, all the input target contraction coefficients are equal to a value $0<\theta ^*\leq 1$ and all the output target dilation coefficients are equal to a value $\phi ^*\geq 1$, (see \eqref{eq:thetalin} and \eqref{eq:thetaquad}). That is,
\begin{equation}
\theta ^*_i=\theta ^*,\quad \phi ^*_r=\phi ^*,\qquad i=1,\ldots ,m,\quad r=1,\ldots ,s.
\end{equation} 
Hence, by \eqref{eq:rho1}, the Farrell oriented efficiency is given by their quotient,
\begin{equation}
\label{eq:rhopc}
\rho ^*=\frac{\frac{1}{m}\sum _{i=1}^m\theta ^*_i}{\frac{1}{s}\sum _{r=1}^s\phi ^*_r}=\frac{\theta ^*}{\phi ^*}.
\end{equation}

All the following results assume CRS. Lemmas \ref{lem1} and \ref{lem2} are general results that do not depend on the particular case we intend to study.

\begin{lemma}
\label{lem1}
Given $\left( \mathbf{x},\mathbf{y}\right) \in \partial ^{\text{W}}P$, we have $\lambda (\mathbf{x},\mathbf{y})\in \partial ^{\text{W}}P$ for all $\lambda >0$.
\end{lemma}

\begin{lemma}
\label{lem2}
Given $\textrm{DMU}_o=(\mathbf{x}_o,\mathbf{y}_o)\in \mathcal{D}$, let $\theta, \theta '\in \left] 0,1\right] $ and $\phi ,\phi '\geq 1$ be such that $\left( \theta \mathbf{x}_o,\phi \mathbf{y}_o\right) ,\left( \theta '\mathbf{x}_o,\phi '\mathbf{y}_o\right) \in \partial ^{\text{W}}P$. Then, $\frac{\theta}{\theta '}=\frac{\phi}{\phi '}$.
\end{lemma}

\begin{proposition}
\label{prop1}
Given $\textrm{DMU}_o=(\mathbf{x}_o,\mathbf{y}_o)\in \mathcal{D}$, let $(\mathbf{x}^*_o,\mathbf{y}^*_o)_{\text{CCR}}$ be the target calculated by the input oriented radial model under CRS (also known as  CCR model) and let $(\mathbf{x}^*_o,\mathbf{y}^*_o)_{\text{G}}$ be the target calculated by a generalized oriented model under CRS with orientation coefficients $d^{\minus}_i=d^{\minus}$ for $i=1,\ldots ,m$, and $d^{\plus}_r=d^{\plus}$ for $r=1,\ldots ,s$, where $d^{\minus},d^{\plus}>0 $. Then, there exists $\lambda \geq 0$ such that $(\mathbf{x}^*_o,\mathbf{y}^*_o)_{\text{G}}=\lambda (\mathbf{x}^*_o,\mathbf{y}^*_o)_{\text{CCR}}$.
\end{proposition}

\begin{corollary}
\label{cor1}
Let $\rho ^*_{\text{G}}$ be the Farrell oriented efficiency given by a generalized oriented model under CRS with orientation coefficients $d^{\minus}_i=d^{\minus}$ for $i=1,\ldots ,m$, and $d^{\plus}_r=d^{\plus}$ for $r=1,\ldots ,s$, where $d^{\minus},d^{\plus}>0$. Then, $\rho ^*_{\text{G}}$ coincides with the classic Farrell efficiency given by the input oriented radial model under CRS (CCR model). 
\end{corollary}

\begin{figure}[tb]
\begin{center}
  \includegraphics[width=0.55\textwidth]{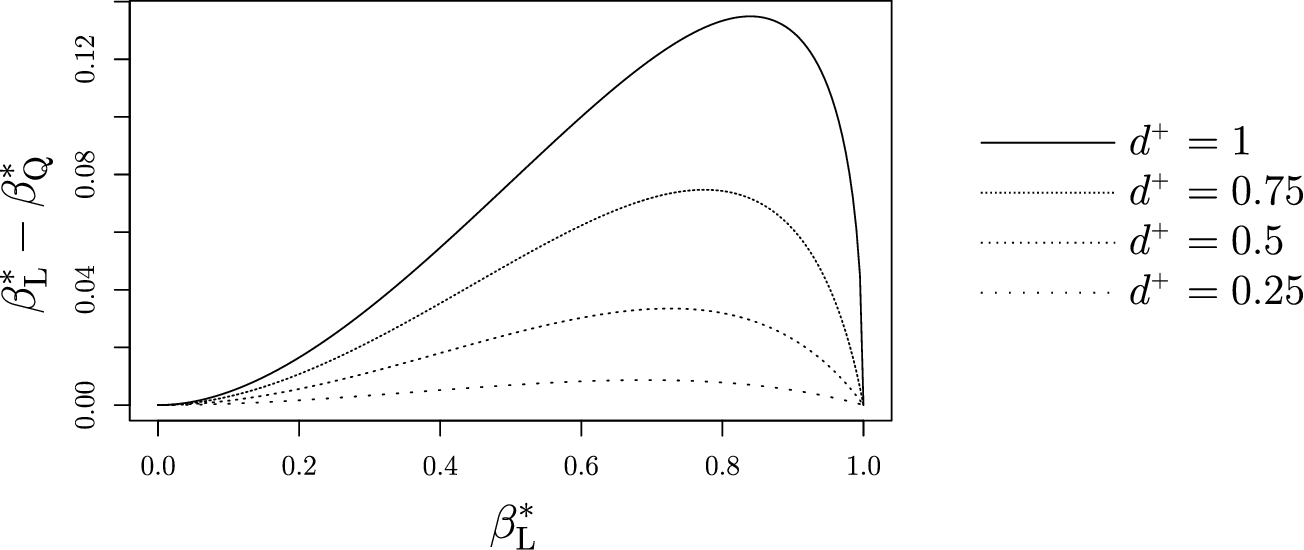}
\end{center}
\caption{Plots of $\beta ^*_{\text{L}}-\beta ^*_{\text{Q}}$ in the particular cases $\mathbf{d}^{\protect \minus}=\mathbf{1}$ and $d^{\protect \plus}_r=d^{\protect \plus}$ for $r=1,\ldots ,s$ under CRS, with different values $d^+=1,0.75,0.5,0.25$, according to \eqref{eq:alphapc}.}
\label{fig:beta1}
\end{figure}

\begin{figure}[tb]
\begin{center}
  \includegraphics[width=0.56\textwidth]{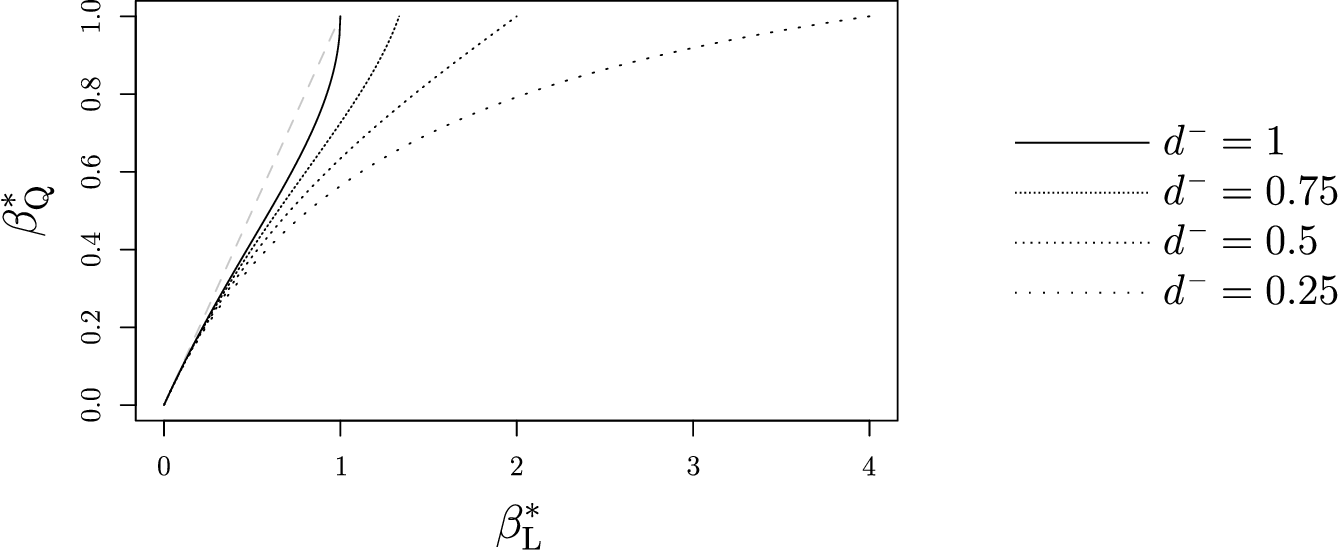}
\end{center}
\caption{Plots of $\beta ^*_{\text{Q}}$ in the particular cases $d^{\protect \minus}_i=d^{\protect \minus}$ for $i=1,\ldots ,m$ and $\mathbf{d}^{\protect \plus}=\mathbf{1}$ under CRS, with different values $d^{\protect \minus}=1,0.75,0.5,0.25$, according to \eqref{eq:alphapc}. Note that $\beta ^*_{\text{L}}$ ranges from $0$ to $1/d^{\protect \minus}$, while $\beta ^*_{\text{Q}}$ ranges from $0$ to $1$.}
\label{fig:beta2}
\end{figure}

Corollary \ref{cor1} assures that, in the important particular case $d^{\minus}_i=d^{\minus}$ for $i=1,\ldots ,m$, and $d^{\plus}_r=d^{\plus}$ for $r=1,\ldots ,s$, with $d^{\minus},d^{\plus}>0$ and under CRS, the Farrell oriented efficiencies $\rho ^*_{\text{L}}$, $\rho ^*_{\text{Q}}$ are equal and they do not depend on the values of $d^{\minus},d^{\plus}$. Moreover, the optimal value $\beta ^*_{\text{Q}}$ and the target $\left( \mathbf{x}^*_o,\mathbf{y}^*_o\right) _{\text{Q}}$ of a QO model can be computed from the optimal value $\beta ^*_{\text{L}}$ of the corresponding LO model, as we are going to show next. From \eqref{eq:thetalin},
\begin{equation}
\label{eq:thetalinpc}
\theta ^*_{\text{L}} = 1-\beta ^*_{\text{L}}d^{\minus},
\qquad \phi ^*_{\text{L}} = 1+\beta ^*_{\text{L}}d^{\plus}.
\end{equation}
On the other hand, from \eqref{eq:thetaquad},
\begin{equation}
\label{eq:thetaquadpc}
\theta ^*_{\text{Q}} = 1-\beta ^*_{\text{Q}}d^{\minus},
\qquad \phi ^*_{\text{Q}} = 1/\left( 1-\beta ^*_{\text{Q}}d^{\plus}\right) .
\end{equation}
Moreover, from Proposition \ref{prop1}, $\frac{\theta ^*_{\text{Q}}}{\theta ^*_{\text{L}}}=\frac{\phi ^*_{\text{Q}}}{\phi ^*_{\text{L}}}$.
Then, from \eqref{eq:thetalinpc} and \eqref{eq:thetaquadpc}, we have $\beta ^{*2}_{\text{Q}}d^{\minus}d^{\plus}-\beta ^*_{\text{Q}}\left( d^{\minus}+d^{\plus}\right) +\beta ^*_{\text{L}}\left( d^{\minus}+d^{\plus}\right) /\left( 1+\beta ^*_{\text{L}}d^{\plus}\right) =0$
and hence,
\begin{equation}
\label{eq:alphapc}
\beta ^*_{\text{Q}} = \frac{d^{\minus}+d^{\plus}-\sqrt{\left( d^{\minus}+d^{\plus}\right) ^2-4\beta ^*_{\text{L}}d^{\minus}d^{\plus}\left( d^{\minus}+d^{\plus}\right) /\left( 1+\beta ^*_{\text{L}}d^{\plus}\right) }}{2d^{\minus}d^{\plus}}.
\end{equation}
From the value of $\beta ^*_{\text{Q}}$ in terms of $\beta ^*_{\text{L}}$ given by \eqref{eq:alphapc} (see Fig. \ref{fig:beta1} and \ref{fig:beta2}), we can find the target $\left( \mathbf{x}^*_o,\mathbf{y}^*_o\right) _{\text{Q}}=\left( \theta ^*_{\text{Q}}\mathbf{x}_o,\phi ^*_{\text{Q}}\mathbf{y}_o\right) $ computing $\theta ^*_{\text{Q}}$ and $\phi ^*_{\text{Q}}$ from \eqref{eq:thetaquadpc}.

In the important particular case of non-oriented radial models, for which $d^{\minus}=d^{\plus}=1$, we have
\[
\beta ^*_{\text{Q}}=1-\sqrt{\frac{1-\beta ^*_{\text{L}}}{1+\beta ^*_{\text{L}}}},\qquad \rho ^*_{\text{Q}}=\rho ^*_{\text{L}}=\frac{1-\beta ^*_{\text{L}}}{1+\beta ^*_{\text{L}}},\qquad
\theta ^*_{\text{Q}}=\sqrt{\frac{1-\beta ^*_{\text{L}}}{1+\beta ^*_{\text{L}}}},\qquad  \phi ^*_{\text{Q}}=\sqrt{\frac{1+\beta ^*_{\text{L}}}{1-\beta ^*_{\text{L}}}}.
\]

\section{Examples}
\label{sec:ex}

In this section, our purpose is to illustrate the defined concepts, the results and the differences between models. Generalized oriented models (as well as directional models) are applicable in any situation where the producer wants to manage the proportions by which inputs and outputs are improved to achieve technical efficiency.

\begin{table}[tb]
\begin{center}
\caption{Values of the variables of the DMUs considered in the example. Only $A$ is efficient under CRS.}
\label{tab1}
\begin{tabular}{c|cc|cc}
  & Input 1 & Input 2 & Output 1 & Output 2 \\ \hline
$A$ & 1       & 1       & 4        & 4        \\
$B$ & 1       & 2       & 1        & 2        \\
$C$ & 1       & 2       & 2        & 1        \\
$D$ & 2       & 1       & 1        & 2        \\
$E$ & 2       & 1       & 2        & 1       
\end{tabular}
\end{center}
\end{table}

Let us consider $\mathcal{D}=\left\{ A,B,C,D,E\right\}$ a set of $5$ DMUs with 2 inputs and 2 outputs given by Table \ref{tab1} and let us assume CRS. In this simple example, only $A$ is efficient under CRS.
The efficient frontier $\partial ^{\text{S}}P$ is formed by all the positive multiples of $A$:
\[
\partial ^{\text{S}}P=\left\{ \alpha \left( 1,1;4,4\right) \,\, |\,\, \alpha >0\right\} .
\]
Moreover, there is always a single efficient projection in $\partial ^{\text{S}}P$ dominating the projection and hence, there is no need for a second stage (see Remark \ref{rem:ss}).

We have used \texttt{R} 4.4.0 \citep{R24} for computations. Specifically, we have used packages \texttt{deaR} \citep{deaR23} for solving LO models, and \texttt{optiSolve} \citep{optisolve} for solving QO models, by means of the augmented Lagrangian minimization (``alabama'') algorithm \citep{lange2004}.

\subsection{Case 1: \texorpdfstring{$\mathbf{d}^{\protect \minus}=\mathbf{d}^{\protect \plus}=(1,1)$}{d-=d+=(1,1)}}

\begin{table}[tb]
\begin{center}
\caption{Results from LO model under CRS with orientation vectors  $\mathbf{d}^{\protect \minus}=\mathbf{d}^{\protect \plus}=(1,1)$.}
\label{tabL1}
{\small
\begin{tabular}{c|cccc}
    & $\beta ^*_{\text{L}}$ & $\rho ^*_{\text{L}}$ & Target in $\partial ^{\text{W}}P$ & Efficient Projection in $\partial ^{\text{S}}P$ \\ \hline
$B=(1,2;1,2)$ & 0.333                 & 0.5                  & $(0.667,1.333;1.333,2.667)$               & $(0.667,0.667;2.667,2.667)$         \\
$C=(1,2;2,1)$ & 0.333                 & 0.5                  & $(0.667,1.333;2.667,1.333)$               & $(0.667,0.667;2.667,2.667)$         \\
$D=(2,1;1,2)$ & 0.333                 & 0.5                  & $(1.333,0.667;1.333,2.667)$               & $(0.667,0.667;2.667,2.667)$         \\
$E=(2,1;2,1)$ & 0.333                 & 0.5                  & $(1.333,0.667;2.667,1.333)$               & $(0.667,0.667;2.667,2.667)$        
\end{tabular}
}
\end{center}
\end{table}

\begin{table}[tb]
\begin{center}
\caption{Results from QO models under CRS with orientation vectors $\mathbf{d}^{\protect \minus}=\mathbf{d}^{\protect \plus}=(1,1)$.}
\label{tabQ1}
{\small
\begin{tabular}{c|cccc}
    & $\beta ^*_{\text{Q}}$ & \hspace{-0.2cm} $\rho ^*_{\text{Q}}$ \hspace{-0.3cm} & Target in $\partial ^{\text{W}}P$ & Efficient Projection in $\partial ^{\text{S}}P$ \\ \hline
$B=(1,2;1,2)$ & 0.293                 & 0.5                  & $(0.707,1.414;1.414,2.828)$               & $(0.707,0.707;2.828,2.828)$         \\
$C=(1,2;2,1)$ & 0.293                 & 0.5                  & $(0.707,1.414;2.828,1.414)$               & $(0.707,0.707;2.828,2.828)$         \\
$D=(2,1;1,2)$ & 0.293                 & 0.5                  & $(1.414,0.707;1.414,2.828)$               & $(0.707,0.707;2.828,2.828)$         \\
$E=(2,1;2,1)$ & 0.293                 & 0.5                  & $(1.414,0.707;2.828,1.414)$               & $(0.707,0.707;2.828,2.828)$        
\end{tabular}
}
\end{center}
\end{table}

Firstly, we consider non-oriented radial models (see Remark \ref{rem:no}), i.e. with orientation vectors $\mathbf{d}^{\minus}=\mathbf{d}^{\plus}=(1,1)$. According to Corollary \ref{cor1}, all models give a Farrell oriented efficiency equal to the classic Farrell efficiency, that is $0.5$ for all inefficient DMUs, as shown in Tables \ref{tabL1} and \ref{tabQ1} (note that values shown in tables are approximated and rounded to three decimal places). It is also noted that $\beta ^*_{\text{Q}}\leq \beta ^*_{\text{L}}$, according to Proposition \ref{propab}.

Although, by symmetry, DMUs $B,C,D,E$ have the same efficient projection in $\partial ^{\text{S}}P$, their targets in $\partial ^{\text{W}}P$ are different. So, according to the LO model, all these DMUs must reduce their inputs to $(2/3,2/3)$ and increase their outputs to $(8/3,8/3)$ to reach the efficient frontier $\partial ^{\text{S}}P$. However, for example, DMU $B$ only has to reduce its inputs to $(2/3,4/3)$ and increase its outputs to $(4/3,8/3)$ for becoming weakly efficient. In this case, any infinitesimal improvement of the first input or the second output (which are variables with inefficiency projection slack equal to zero, see Remark \ref{rem:ss}) would transform this weakly efficient activity into an efficient one. On the other hand, according to the QO model, all the DMUs must reduce their inputs to $(0.707,0.707)$ approx., and increase their outputs to $(2.828,2.828)$ approx., to reach the efficient frontier, and for example, DMU $B$ only has to reduce its inputs to $(0.707,1.414)$ approx., and increase its outputs to $(1.414,2.828)$ approx., for becoming weakly efficient. 

With respect to the balance between input contractions and output dilations in the calculation of the target in $\partial ^{\text{W}}P$, the improvements proposed by the LO model are not according to CRS (see Remark \ref{rem:crs}). For example, DMU $B$ (see Table \ref{tabB}) must contract its inputs by $\theta ^*_{\text{L}}\approx 0.667$ and dilate its outputs by $\phi ^*_{\text{L}}\approx 1.333$ for becoming weakly efficient, while the dilation according to CRS would be $1.5$ (i.e. the inverse of $\theta ^*_{\text{L}}$). In contrast, the proposed relative slacks are the same for inputs and for outputs, $\tau^{\minus *}_{\text{L}}=\tau ^{\plus *}_{\text{L}}\approx 0.333$.

On the other hand, the QO model proposes different relative slacks for inputs, $\tau ^{\minus *}_{\text{Q}}\approx 0.293$, and for outputs, $\tau ^{\plus *}_{\text{Q}}\approx 0.414$. However, this model proposes a contraction of inputs $\theta ^*_{\text{Q}}\approx 0.707$ and a dilation of outputs given by its inverse, i.e. $\phi ^*_{\text{Q}}\approx 1.414$, according to CRS.

\subsection{Case 2: \texorpdfstring{$\mathbf{d}^{\protect \minus}=(1,1)$ and $\mathbf{d}^{\protect \plus}=(0.5,0.5)$}{d-=(1,1) and d+=(0.5,0.5)}}

\begin{table}[tb]
\begin{center}
\caption{Results from LO model under CRS with orientation vectors $\mathbf{d}^{\protect \minus}=(1,1)$ and $\mathbf{d}^{\protect \plus}=(0.5,0.5)$.}
\label{tabL2}
{\small
\begin{tabular}{c|cccc}
    & $\beta ^*_{\text{L}}$ & $\rho ^*_{\text{L}}$ & Target in $\partial ^{\text{W}}P$ & Efficient Projection in $\partial ^{\text{S}}P$ \\ \hline
$B=(1,2;1,2)$ & 0.4                 & 0.5                  & $(0.6,1.2;1.2,2.4)$               & $(0.6,0.6;2.4,2.4)$         \\
$C=(1,2;2,1)$ & 0.4                 & 0.5                  & $(0.6,1.2;2.4,1.2)$               & $(0.6,0.6;2.4,2.4)$         \\
$D=(2,1;1,2)$ & 0.4                 & 0.5                  & $(1.2,0.6;1.2,2.4)$               & $(0.6,0.6;2.4,2.4)$         \\
$E=(2,1;2,1)$ & 0.4                 & 0.5                  & $(1.2,0.6;2.4,1.2)$               & $(0.6,0.6;2.4,2.4)$        
\end{tabular}
}
\end{center}
\end{table}

\begin{table}[tb]
\begin{center}
\caption{Results from QO model under CRS with orientation vectors $\mathbf{d}^{\protect \minus}=(1,1)$ and $\mathbf{d}^{\protect \plus}=(0.5,0.5)$.}
\label{tabQ2}
{\small
\begin{tabular}{c|cccc}
    & $\beta ^*_{\text{Q}}$ & $\rho ^*_{\text{Q}}$ & Target in $\partial ^{\text{W}}P$ & Efficient Projection in $\partial ^{\text{S}}P$ \\ \hline
$B=(1,2;1,2)$ & 0.382                 & 0.5                  & $(0.618,1.236;1.236,2.472)$               & $(0.618,0.618;2.472,2.472)$         \\
$C=(1,2;2,1)$ & 0.382                 & 0.5                  & $(0.618,1.236;2.472,1.236)$               & $(0.618,0.618;2.472,2.472)$         \\
$D=(2,1;1,2)$ & 0.382                 & 0.5                  & $(1.236,0.618;1.236,2.472)$               & $(0.618,0.618;2.472,2.472)$         \\
$E=(2,1;2,1)$ & 0.382                 & 0.5                  & $(1.236,0.618;2.472,1.236)$               & $(0.618,0.618;2.472,2.472)$        
\end{tabular}
}
\end{center}
\end{table}

Secondly, we consider orientation vectors $\mathbf{d}^{\minus}=(1,1)$ and $\mathbf{d}^{\plus}=(0.5,0.5)$. As in the first case and according to Corollary \ref{cor1}, all models give the same Farrell oriented efficiency, equal to the classic Farrell efficiency $0.5$, as shown in Tables \ref{tabL2} and \ref{tabQ2}. Moreover, by symmetry, the values of parameter $\beta ^*$ do not depend on the evaluated DMU ($B,C,D,E$) but on the model, satisfying $\beta ^*_{\text{Q}}\leq \beta ^*_{\text{L}}$, according to Proposition \ref{propab}. We note that DMUs $B,C,D,E$ have the same efficient projection in $\partial ^{\text{S}}P$, but their targets in $\partial ^{\text{W}}P$ are different, as in the first case.

With respect to the balance between input contractions and output dilations in the calculation of the target in $\partial ^{\text{W}}P$, the LO model proposes for DMU $B$ (see Table \ref{tabB}) a contraction of inputs $\theta ^*_{\text{L}}=0.6$ and a dilation of outputs $\phi ^*_{\text{L}}=1.2$, or equivalently, proposes improving inputs applying a relative slack of $\tau ^{\minus *}_{\text{L}}=0.4$ and improving outputs applying a relative slack of $\tau ^{\plus *}_{\text{L}}=0.2$ (i.e. $0.5$ times $\tau ^{\minus *}_{\text{L}}$, where $0.5$ corresponds to the output orientation coefficient). 

For its part, the QO model proposes a contraction of inputs $\theta ^*_{\text{Q}}\approx 0.618$ and a dilation of outputs $\phi ^*_{\text{Q}}\approx 1.236$, that is the inverse of the contraction that would have an input with orientation coefficient equal to $0.5$, i.e. the inverse of $1-(1-\theta ^*_{\text{Q}})\cdot 0.5$.

\subsection{Case 3: \texorpdfstring{$\mathbf{d}^{\protect \minus}=\mathbf{d}^{\protect \plus}=(1,0.5)$}{d-=d+=(1,0.5)}}

\begin{table}[tb]
\begin{center}
\caption{Results from LO model under CRS with orientation vectors $\mathbf{d}^{\protect \minus}=\mathbf{d}^{\protect \plus}=(1,0.5)$.}
\label{tabL3}
{\small
\begin{tabular}{c|cccc}
    & $\beta ^*_{\text{L}}$ & $\rho ^*_{\text{L}}$ & Target in $\partial ^{\text{W}}P$ & Efficient Projection in $\partial ^{\text{S}}P$ \\ \hline
$B=(1,2;1,2)$ & 0.4                 & 0.538                  & $(0.6,1.6;1.4,2.4)$               & $(0.6,0.6;2.4,2.4)$         \\
$C=(1,2;2,1)$ & 0.333                 & 0.6                  & $(0.667,1.667;2.667,1.167)$               & $(0.667,0.667;2.667,2.667)$         \\
$D=(2,1;1,2)$ & 0.667                 & 0.333                  & $(0.667,0.667;1.667,2.667)$               & $(0.667,0.667;2.667,2.667)$         \\
$E=(2,1;2,1)$ & 0.5                 & 0.455                  & $(1,0.75;3,1.25)$               & $(0.75,0.75;3,3)$        
\end{tabular}
}
\end{center}
\end{table}

\begin{table}[tb]
\begin{center}
\caption{Results from QO model under CRS with orientation vectors $\mathbf{d}^{\protect \minus}=\mathbf{d}^{\protect \plus}=(1,0.5)$.}
\label{tabQ3}
{\small
\begin{tabular}{c|cccc}
    & $\beta ^*_{\text{Q}}$ & $\rho ^*_{\text{Q}}$ & Target in $\partial ^{\text{W}}P$ & Efficient Projection in $\partial ^{\text{S}}P$ \\ \hline
$B=(1,2;1,2)$ & 0.382                 & 0.5                  & $(0.618,1.618;1.618,2.472)$               & $(0.618,0.618;2.472,2.472)$         \\
$C=(1,2;2,1)$ & 0.293                 & 0.604                  & $(0.707,1.707;2.828,1.172)$               & $(0.707,0.707;2.828,2.828)$         \\
$D=(2,1;1,2)$ & 0.586                 & 0.293                  & $(0.828,0.707;2.414,2.828)$               & $(0.707,0.707;2.828,2.828)$         \\
$E=(2,1;2,1)$ & 0.382                 & 0.5                  & $(1.236,0.809;3.236,1.236)$               & $(0.809,0.809;3.236,3.236)$        
\end{tabular}
}
\end{center}
\end{table}

Thirdly, we consider orientation vectors $\mathbf{d}^{\minus}=\mathbf{d}^{\plus}=(1,0.5)$. In this case, as shown in Tables \ref{tabL3} and \ref{tabQ3}, symmetry is broken and the values of $\beta ^*$ and $\rho ^*$ depend on the model but also on the evaluated DMU.
According to Proposition \ref{propab}, $\beta ^*_{\text{Q}}\leq \beta ^*_{\text{L}}$ holds for any DMU.
Moreover, contrary to the previous cases, not all DMUs have the same efficient projection in $\partial ^{\text{S}}P$, although some DMUs do have the same, as in the case of $C$ and $D$ in LO and QO models. However, all their targets in $\partial ^{\text{W}}P$ are different.

For this orientation, $C$ has the highest Farrell oriented efficiency in all models, followed by $B$, $E$ and $D$. Note that in the QO model, there is a tie between $B$ and $E$, both with $\rho ^*_{\text{Q}}=0.5$. 

\begin{table}[tb]
\begin{center}
\caption{Input target contractions and output target dilations vectors, $\bm{\theta}^*,\bm{\phi}^*$, and relative input and output target slacks vectors, $\bm{\tau}^{{\protect \minus}*},\bm{\tau}^{{\protect \plus}*}$ of DMU $B=(1,2;1,2)$ evaluated by LO and QO models under CRS with different orientation vectors $\mathbf{d}^{\protect \minus},\mathbf{d}^{\protect \plus}$.}
\label{tabB}
{\small
\begin{tabular}{r|cc|cccc}
        & \hspace{-0.2cm} $\mathbf{d}^{\protect \minus}$             & \hspace{-0.4cm} $\mathbf{d}^{\protect \plus}$ & $\bm{\theta}^*$ & \hspace{-0.2cm} $\bm{\phi}^*$ & \hspace{-0.2cm} $\bm{\tau}^{{\protect \minus}*}$ & \hspace{-0.2cm} $\bm{\tau}^{{\protect \plus}*}$ \\ \hline
LO      & \hspace{-0.2cm} \multirow{2}{*}{$(1,1)$}   & \hspace{-0.4cm} \multirow{2}{*}{$(1,1)$} & $(0.667,0.667)$ & \hspace{-0.2cm} $(1.333,1.333)$ & \hspace{-0.2cm} $(0.333,0.333)$ & \hspace{-0.2cm} $(0.333,0.333)$  \\
QO &                            &                              & $(0.707,0.707)$ & \hspace{-0.2cm} $(1.414,1.414)$ & \hspace{-0.2cm} $(0.293,0.293)$  & \hspace{-0.2cm} $(0.414,0.414)$  \\ \hline
LO      & \hspace{-0.2cm} \multirow{2}{*}{$(1,1)$}   & \hspace{-0.4cm} \multirow{2}{*}{$(0.5,0.5)$} & $(0.6,0.6)$ & \hspace{-0.2cm} $(1.2,1.2)$ & \hspace{-0.2cm} $(0.4,0.4)$ & \hspace{-0.2cm} $(0.2,0.2)$      \\
QO      &                            &                              & $(0.618,0.618)$ & \hspace{-0.2cm} $(1.236,1.236)$ & \hspace{-0.2cm} $(0.382,0.382)$  & \hspace{-0.2cm} $(0.236,0.236)$  \\ \hline
LO      & \hspace{-0.2cm} \multirow{2}{*}{$(1,0.5)$} & \hspace{-0.4cm} \multirow{2}{*}{$(1,0.5)$}   & $(0.6,0.8)$     & \hspace{-0.2cm} $(1.4,1.2)$     & \hspace{-0.2cm} $(0.4,0.2)$      & \hspace{-0.2cm} $(0.4,0.2)$      \\
QO      &                            &                              & $(0.618,0.809)$ & \hspace{-0.2cm} $(1.618,1.236)$ & \hspace{-0.2cm} $(0.382,0.191)$  & \hspace{-0.2cm} $(0.618,0.236)$ 
\end{tabular}
}
\end{center}
\end{table}

With respect to the balance between input contractions and output dilations in the calculation of the projection onto $\partial ^{\text{W}}P$, the LO model proposes improving the variables by means of relative slacks vectors which are directly proportional to the orientation vectors with the same proportionality factor $\beta ^*_{\text{L}}$. For example, for DMU $B$ (see Table \ref{tabB}), the relative slacks vectors are given by $\bm{\tau}^{\minus *}_{\text{L}}=0.4\cdot \mathbf{d}^{\minus}=(0.4,0.2)$ and $\bm{\tau}^{\plus *}_{\text{L}}=0.4\cdot \mathbf{d}^{\plus}=(0.4,0.2)$. As a consequence, the dilation of the first output, $\phi ^*_{\text{L}\, 1}=1.4$, is not the inverse of the contraction of the first input, $\theta ^*_{\text{L}\, 1}=0.6$, both variables with orientation coefficient equal to $1$, showing that this model is not balanced according to CRS.

On the other hand, the QO model proposes output dilations following the CRS balance and, moreover, the dilation of an output with orientation coefficient $d$ is the inverse of the contraction that would have an input with the same orientation coefficient $d$. For example, for DMU $B$ (see Table \ref{tabB}), the dilation of the first output, $\phi ^*_{\text{Q}\, 1}\approx 1.618$ (golden ratio, by the way), is the inverse of the contraction of the first input, $\theta ^*_{\text{Q}\, 1}\approx 0.618$, both variables with orientation coefficient equal to $1$; and the dilation of the second output, $\phi ^*_{\text{Q}\, 2}\approx 1.236$, is the inverse of the contraction of the second input, $\theta ^*_{\text{Q}\, 2}\approx 0.809$, both variables with orientation coefficient equal to $0.5$.

\section{Final remarks}
\label{sec:fin}

We have reviewed linear directional (LO) models and introduced QO models. All of these models employ a generalized orientation $\left( \mathbf{d}^{\minus},\mathbf{d}^{\plus}\right) $ for calculating a target in the weakly efficient frontier of the production possibility set. However, each model is applicable in specific circumstances. LO models have the particularity that the relative target slacks vectors (both, input and output) are always directly proportional to the orientation vectors with the same proportionality factor, $\beta ^*_{\text{L}}$ (see \eqref{eq:taulin}). Due to their linearity, LO models are appropriate when either $\mathbf{d}^{\minus}=\mathbf{0}$ or $\mathbf{d}^{\plus}=\mathbf{0}$, indicating that there is no desire to simultaneously improve inputs and outputs.
Otherwise, it is necessary to consider that LO models are not balanced according to CRS (see Remark \ref{rem:crs}) and thus they are not suitable for simultaneously improving inputs and outputs under the CRS assumption.
In these cases, we have demonstrated that QO models should be selected because they are balanced according to CRS. In fact, for QO models, the target dilation of an output with orientation coefficient $d$ is the inverse of the target contraction that would have an input with equal orientation coefficient $d$ (see \eqref{eq:thetaquad}).

Furthermore, we have broadened the conventional Farrell efficiency to encompass all these generalized oriented models. In particular, we have defined the \emph{Farrell oriented efficiency} as a score that measures technical efficiency  in accordance with the improvement strategy determined by the orientation. Each generalized oriented model has its own Farrell oriented efficiency but, in the particular cases in which $d^{\minus}_i=d^{\minus}$ for $i=1,\ldots ,m$, and $d^{\plus}_r=d^{\plus}$ for $r=1,\ldots ,s$ (with at least one non-zero parameter $d^-$ or $d^+$) and under CRS, they are equal to the classic Farrell efficiency given by the CCR model (see Corollary \ref{cor1}). These particular cases are those in which there is no discrimination between inputs nor between outputs; that is to say, the same contraction is applied to all inputs, and the same dilation is applied to all outputs.

In terms of computational aspects, LO models are linear, whereas QO models are quadratically constrained with a convex feasible set, which can be expressed as convex optimization problems in abstract form. However, they are not second-order cone programs (SOCP). This is not currently a significant disadvantage, as there are algorithms that are highly effective in addressing these non-linear problems, such as the augmented Lagrangian minimization algorithm \citep{lange2004}. In our example, the computational cost of QO models using this algorithm is usually $20$ times greater than that of LO models, although this may vary depending on the number of DMUs and the number of variables.
Nevertheless, we have demonstrated that, in the above particular cases ($d^{\minus}_i=d^{\minus}$ for $i=1,\ldots ,m$, and $d^{\plus}_r=d^{\plus}$ for $r=1,\ldots ,s$), QO models can be solved by means of the LO model, that is to say, in a linear manner.

Future research may consist of adapting generalized oriented models to a stochastic scenario, following the methodology of \citet{Bolos2024} and defining chance-constrained versions of generalized oriented models.
Note that chance-constrained versions of QO models will also be quadratically constrained with convex feasible set and therefore, these models will not lead to a significant increase in computational complexity with respect to deterministic QO models, or with respect to chance constrained LO models.
On the other hand, non-convex technologies are another interesting topic for future research. For this purpose, the methodology given by \citet{Tone2014} about non-convex meta-frontiers (resulting from the categorization of DMUs into several classes) can be adapted to generalized oriented models. This methodology uses some models (radial or non-radial) for computing \textit{scale-efficiencies} as the quotient of CRS and VRS efficiencies, and determine some \textit{scale-dependent} and \textit{scale-independent slacks}.
In this way, Farrell oriented efficiencies from QO and LO models could be used for computing these scale efficiencies, thus taking into account the criteria of management chosen by the producer.


\bibliography{refs}
\bibliographystyle{abbrvnat}


\section*{Appendix A. Proofs}
\label{sec:proofs}

\begin{proof}[\textbf{Proof of Proposition \ref{propab}}]
Reasoning by \textit{reductio ad absurdum}, let us suppose that $\beta ^*_{\text{Q}}>\beta ^*_{\text{L}}$. Then, taking into account that
\begin{equation}
\label{eq:proof31}
1+\beta ^*_{\text{Q}}d^{\plus}_r\leq \frac{1}{1-\beta ^*_{\text{Q}}d^{\plus}_r},\quad\textrm{ for }\quad 0\leq \beta ^*_{\text{Q}}d^{\plus}_r<1,
\end{equation}
we have $\left( \mathrm{diag}\left( 1-\beta ^*_{\text{Q}}\mathbf{d}^{\minus}\right) \mathbf{x}_o,\mathrm{diag}\left( 1+\beta ^*_{\text{Q}}\mathbf{d}^{\plus}\right) \mathbf{y}_o\right) \in P$,
because it is dominated by the target of the QO model. Hence, $\beta ^*_{\text{Q}}$ satisfies the constraints of \eqref{eq:dirgenlin}, leading to $\beta ^*_{\text{Q}}\leq \beta ^*_{\text{L}}$ and a contradiction.
Moreover, since the equality in \eqref{eq:proof31} holds for any possible value of $\beta ^*_{\text{Q}}$ if and only if $d^{\plus}_r=0$, we have that $\beta ^*_{\text{Q}}=\beta ^*_{\text{L}}$ for any evaluated DMU if and only if $\mathbf{d}^{\plus}=\mathbf{0}$.
\end{proof}

\begin{proof}[\textbf{Proof of Proposition \ref{prop:impr}}]
We are going to suppose that all the variables are controllable.
For LO models, taking into account \eqref{eq:imprLQO} and \eqref{eq:taulin}, we have that
\[
f\left( \bm{\tau}^{\minus *}_{\mathrm{L}},\bm{\tau}^{\plus *}_{\mathrm{L}}\right) \approx \nabla f\left( \mathbf{0},\mathbf{0}\right) \cdot \left( \bm{\tau}^{\minus *}_{\mathrm{L}},\bm{\tau}^{\plus *}_{\mathrm{L}}\right) = \nabla f\left( \mathbf{0},\mathbf{0}\right) \cdot \left( \beta ^*_{\mathrm{L}}\mathbf{d}^{\minus},\beta ^*_{\mathrm{L}}\mathbf{d}^{\plus}\right) =\beta ^*_{\mathrm{L}}
,
\]
where $\bm{\tau}^{\minus *}_{\mathrm{L}},\bm{\tau}^{\plus *}_{\mathrm{L}}$ are the relative target slacks vectors.

For QO models, we consider the function $g\left( \bm{\tau}^{\minus},\mathbf{b}\right) =f\left( \bm{\tau}^{\minus}, \frac{\mathbf{b}}{1-\mathbf{b}}\right) $, where $0\leq \tau ^{\minus}_i,b_r<1$ for $i=1,\ldots ,m$ and $r=1,\ldots ,s$. Taking into account \eqref{eq:imprLQO}, \eqref{eq:tauquad} and the fact that $\nabla g\left( \mathbf{0},\mathbf{0}\right) =\nabla f\left( \mathbf{0},\mathbf{0}\right) $, we have that
\begin{equation*}
\begin{array}{rl}
f\left( \bm{\tau}^{\minus *}_{\mathrm{Q}},\bm{\tau}^{\plus *}_{\mathrm{Q}}\right) & = g\left( \bm{\tau}^{\minus *}_{\mathrm{Q}},\frac{\bm{\tau}^{\plus *}_{\mathrm{Q}}}{1+\bm{\tau}^{\plus *}_{\mathrm{Q}}}\right) 
\approx \nabla g\left( \mathbf{0},\mathbf{0}\right) \cdot \left( \bm{\tau}^{\minus *}_{\mathrm{Q}},\frac{\bm{\tau}^{\plus *}_{\mathrm{Q}}}{1+\bm{\tau}^{\plus *}_{\mathrm{Q}}}\right) \\  \\
& = \nabla f\left( \mathbf{0},\mathbf{0}\right) \cdot \left( \beta ^*_{\mathrm{Q}}\mathbf{d}^{\minus},\beta ^*_{\mathrm{Q}}\mathbf{d}^{\plus}\right) = \beta ^*_{\mathrm{Q}}
,
\end{array}
\end{equation*}
where $\bm{\tau}^{\minus *}_{\mathrm{Q}},\bm{\tau}^{\plus *}_{\mathrm{Q}}$ are the relative target slacks vectors.
\end{proof}

\begin{proof}[\textbf{Proof of Proposition \ref{prop:monbeta}}]
Let us consider $\beta ^*_{\text{L}}(\mathbf{x},\mathbf{y})$ given by \eqref{eq:dirmodlinxy}. In this program, we replace the value of the first input $x_1$ by a new variable $x$ and add a new constraint $x=x_1$. Hence, the Karush-Kuhn-Tucker (KKT) multiplier of this new constraint is $\frac{\partial \beta ^*_{\text{L}}}{\partial x_1}(\mathbf{x},\mathbf{y})$, and the stationarity condition associated to the new variable $x$ is given by
\[
(1-\beta d^{\minus}_1)\mu _1-\frac{\partial \beta ^*_{\text{L}}}{\partial x_1}(\mathbf{x},\mathbf{y})=0,
\]
where $\mu _1\geq 0$ is the KKT multiplier of the first constraint. Since $0\leq \beta ^*_{\text{L}}\left\| \mathbf{d}^{\minus}\right\| _{\infty}< 1$, we can suppose $\beta d^{\minus}_1<1$, concluding that
\[
\frac{\partial \beta ^*_{\text{L}}}{\partial x_1}(\mathbf{x},\mathbf{y})=(1-\beta d^{\minus}_1)\mu _1\geq 0.
\]

In the same way, the result can be proved for the rest of the inputs.

With respect to the outputs, we replace the value of the first output $y_1$ by a new variable $y$ and add a new constraint $y=y_1$. Hence, the KKT multiplier of this new constraint is $\frac{\partial \beta ^*_{\text{L}}}{\partial y_1}(\mathbf{x},\mathbf{y})$, and the stationarity condition associated to the new variable $y$ is given by
\[
(1+\beta d^{\plus}_1)\nu _1-\frac{\partial \beta ^*_{\text{L}}}{\partial y_1}(\mathbf{x},\mathbf{y})=0,
\]
where $\nu _1\leq 0$ is the KKT multiplier of the $(m+1)$-th constraint (i.e. the first constraint of the block of outputs). Since $\beta \geq 0$, we conclude that 
\[
\frac{\partial \beta ^*_{\text{L}}}{\partial y_1}(\mathbf{x},\mathbf{y})=(1+\beta d^{\plus}_1)\nu _1\leq 0.
\]
In the same way, the result can be proved for the rest of the outputs.

Using the same technique, the result can be proved for $\beta ^*_{\text{Q}}(\mathbf{x},\mathbf{y})$ given by \eqref{eq:dirmodquadxy}. We obtain
\[
\frac{\partial \beta ^*_{\text{Q}}}{\partial x_1}(\mathbf{x},\mathbf{y})=(1-\beta d^{\minus}_1)\mu _1\geq 0
\]
and
\[
\frac{\partial \beta ^*_{\text{Q}}}{\partial y_1}(\mathbf{x},\mathbf{y})=\left( \frac{1}{1-\beta d^{\plus}_1}\right) \nu _1\leq 0,
\]
where $\mu _1\geq 0$ and $\nu _1\leq 0$ are the multipliers of the first and the $(m+1)$-th constraints of \eqref{eq:dirmodquadxy}, respectively.  In the same way, the result can be proved for the rest of variables.
\end{proof}

\begin{proof}[\textbf{Proof of Proposition \ref{prop:monrho}}]
In this proof, we are going to use the following notation:
\[
\overline{d^{\minus}}=\frac{1}{m}\sum _{i=1}^m d^{\minus}_i\qquad \textrm{and}\qquad \overline{d^{\plus}}=\frac{1}{s}\sum _{j=1}^s d^{\plus}_j.
\]
Moreover, we are going to suppose without loss of generality that $\| \left( \mathbf{d}^{\minus},\mathbf{d}^{\plus}\right) \| _{\infty}=1$. 
Let us consider $\rho ^*_{\text{L}}(\mathbf{x},\mathbf{y})$ given by \eqref{eq:rhoxy} with $\beta ^*(\mathbf{x},\mathbf{y})=\beta ^*_{\text{L}}(\mathbf{x},\mathbf{y})$. From \eqref{eq:rholin}, we have
\[
\rho ^*_{\text{L}}(\mathbf{x},\mathbf{y}) =\frac{1-\beta ^*_{\text{L}}(\mathbf{x},\mathbf{y})\overline{d^{\minus}}}{1+\beta ^*_{\text{L}}(\mathbf{x},\mathbf{y})\overline{d^{\plus}}},
\]
and hence,
\[
\frac{\partial \rho ^*_{\text{L}}}{\partial x_1}(\mathbf{x},\mathbf{y})=
-\frac{\overline{d^{\minus}}+\overline{d^{\plus}}}{\left( 1+\beta ^*_{\text{L}}(\mathbf{x},\mathbf{y})\overline{d^{\plus}}\right) ^2}\frac{\partial \beta ^*_{\text{L}}}{\partial x_1}(\mathbf{x},\mathbf{y})
\leq 0,
\]
since $\frac{\partial \beta ^*_{\text{L}}}{\partial x_1}(\mathbf{x},\mathbf{y})\geq 0$ from Proposition \ref{prop:monbeta}. In the same way, the result can be proved for the rest of the inputs.

With respect to the outputs, we have
\[
\frac{\partial \rho ^*_{\text{L}}}{\partial y_1}(\mathbf{x},\mathbf{y})=
-\frac{\overline{d^{\minus}}+\overline{d^{\plus}}}{\left( 1+\beta ^*_{\text{L}}(\mathbf{x},\mathbf{y})\overline{d^{\plus}}\right) ^2}\frac{\partial \beta ^*_{\text{L}}}{\partial y_1}(\mathbf{x},\mathbf{y})
\geq 0,
\]
since $\frac{\partial \beta ^*_{\text{L}}}{\partial y_1}(\mathbf{x},\mathbf{y})\leq 0$ from Proposition \ref{prop:monbeta}. In the same way, the result can be proved for the rest of the outputs.

Let us consider $\rho ^*_{\text{Q}}(\mathbf{x},\mathbf{y})$ given by \eqref{eq:rhoxy} with $\beta ^*(\mathbf{x},\mathbf{y})=\beta ^*_{\text{Q}}(\mathbf{x},\mathbf{y})$. From \eqref{eq:rhoquad}, we have
\[
\rho ^*_{\text{Q}}(\mathbf{x},\mathbf{y}) =
\frac{1-\beta ^*_{\text{Q}}(\mathbf{x},\mathbf{y})\overline{d^{\minus}}}{\displaystyle \frac{1}{s}\sum _{r=1}^s \frac{1}{1-\beta ^*_{\text{Q}}(\mathbf{x},\mathbf{y})d^{\plus}_r}},
\]
and hence
\[
\frac{\partial \rho ^*_{\text{Q}}}{\partial x_1}(\mathbf{x},\mathbf{y})=
-\frac{\displaystyle \sum _{r=1}^s\frac{1-\left( 2\beta ^*_{\text{Q}}(\mathbf{x},\mathbf{y}) -1\right) d^{\plus}_r}{\left( 1-\beta ^*_{\text{Q}}(\mathbf{x},\mathbf{y}) d^{\plus}_r\right) ^2}}
{\displaystyle \left( \sum _{r=1}^s \frac{1}{1-\beta ^*_{\text{Q}}(\mathbf{x},\mathbf{y}) d^{\plus}_r}\right) ^2}
s\overline{d^{\minus}}
\frac{\partial \beta ^*_{\text{Q}}}{\partial x_1}
(\mathbf{x},\mathbf{y})
\leq 0,
\]
since $\frac{\partial \beta ^*_{\text{Q}}}{\partial x_1}(\mathbf{x},\mathbf{y})\geq 0$ from Proposition \ref{prop:monbeta}, and $\beta ^*_{\text{Q}}(\mathbf{x},\mathbf{y})<1$. In the same way, the result can be proved for the rest of the inputs.

With respect to the outputs, we have
\[
\frac{\partial \rho ^*_{\text{Q}}}{\partial y_1}(\mathbf{x},\mathbf{y})=
-\frac{\displaystyle \sum _{r=1}^s\frac{1-\left( 2\beta ^*_{\text{Q}}(\mathbf{x},\mathbf{y}) -1\right) d^{\plus}_r}{\left( 1-\beta ^*_{\text{Q}}(\mathbf{x},\mathbf{y}) d^{\plus}_r\right) ^2}}
{\displaystyle \left( \sum _{r=1}^s \frac{1}{1-\beta ^*_{\text{Q}}(\mathbf{x},\mathbf{y}) d^{\plus}_r}\right) ^2}
s\overline{d^{\minus}}
\frac{\partial \beta ^*_{\text{Q}}}{\partial y_1}
(\mathbf{x},\mathbf{y})
\geq 0,
\]
since $\frac{\partial \beta ^*_{\text{Q}}}{\partial y_1}(\mathbf{x},\mathbf{y})\leq 0$ from Proposition \ref{prop:monbeta}, and $\beta ^*_{\text{Q}}(\mathbf{x},\mathbf{y})<1$. In the same way, the result can be proved for the rest of the outputs.

\end{proof}

\begin{proof}[\textbf{Proof of Lemma \ref{lem1}}]
Reasoning by \textit{reductio ad absurdum}, let us suppose that there exists $\lambda >0$ such that $\lambda (\mathbf{x},\mathbf{y})\notin \partial ^{\text{W}}P$. According to \eqref{eq:partialW}, there exists $(\mathbf{x}',\mathbf{y}')\in P$ such that $x'_i<\lambda x_i$ and $y'_r>\lambda y_r$, for $i=1,\ldots ,m$ and $r=1,\ldots ,s$. Hence, the activity $\frac{1}{\lambda} (\mathbf{x}',\mathbf{y}')$ is in $P$ (under CRS) and improves all the variables of $(\mathbf{x},\mathbf{y})$, because $x'_i/\lambda <x_i$ and $y'_r/\lambda >y_r$ for $i=1,\ldots ,m$ and $r=1,\ldots ,s$, leading to a contradiction with $(\mathbf{x},\mathbf{y})\in \partial ^{\text{W}}P$. 
\end{proof}

\begin{proof}[\textbf{Proof of Lemma \ref{lem2}}]
Let us define $\lambda \defeq \theta /\theta '$.
We are going to prove that $\phi =\lambda \phi '$.
Reasoning by \textit{reductio ad absurdum}, let us suppose that $\phi >\lambda \phi '$. Taking $\epsilon >0$ sufficiently small such that $(1-\epsilon)\phi >\lambda \phi '$, the activity given by $\left( 1-\epsilon\right) \left( \theta \mathbf{x}_o,\phi \mathbf{y}_o\right) $
is in $P$ (under CRS, in fact it is in $\partial ^{\text{W}}P$ from Lemma \ref{lem1}) and improves all the variables of $\left( \theta \mathbf{x}_o,\lambda \phi '\mathbf{y}_o\right) $. But, from the definition of $\lambda$ and Lemma \ref{lem1}, we have
\[
\left( \theta \mathbf{x}_o,\lambda \phi '\mathbf{y}_o\right) =\lambda \left( \theta '\mathbf{x}_o,\phi '\mathbf{y}_o\right) \in \partial ^{\text{W}}P,
\]
leading to a contradiction.

On the other hand, let us suppose that $\phi <\lambda \phi '$. Taking $\epsilon >0$ sufficiently small such that $\phi <(1-\epsilon )\lambda \phi '$, the activity given by $\left( 1-\epsilon \right) \left( \theta \mathbf{x}_o,\lambda \phi '\mathbf{y}_o\right) =\left( 1-\epsilon \right) \lambda \left( \theta '\mathbf{x}_o,\phi '\mathbf{y}_o\right) $
is in $P$ (under CRS, in fact, it is in $\partial ^{\text{W}}P$ from Lemma \ref{lem1}) and improves all the variables of $\left( \theta \mathbf{x}_o,\phi \mathbf{y}_o\right) \in \partial ^{\text{W}}P$, leading to a contradiction.
\end{proof}

\begin{proof}[\textbf{Proof of Proposition \ref{prop1}}]
Let us consider the activities $(\mathbf{x}^*_o,\mathbf{y}^*_o)_\text{G}=(\theta ^*_{\text{G}}\mathbf{x}_o,\phi ^*_{\text{G}}\mathbf{y}_o)$ and $(\mathbf{x}^*_o,\mathbf{y}^*_o)_\text{CCR}=(\theta ^*_{\text{CCR}}\mathbf{x}_o,\mathbf{y}_o)$, both in $\partial ^{\text{W}}P$. From Lemma \ref{lem2}, we take $\lambda = \theta ^*_{\text{G}}/\theta ^*_{\text{CCR}}=\phi ^*_{\text{G}}$.
\end{proof}

\begin{proof}[\textbf{Proof of Corollary \ref{cor1}}]
From \eqref{eq:rhopc} and Proposition \ref{prop1}, we have $\rho ^*_{\text{G}}=\theta ^*_{\text{G}}/\phi ^*_{\text{G}}=\theta ^*_{\text{CCR}}$.
\end{proof}

\end{document}